\documentclass{article}
\usepackage{latexsym}
\usepackage{amsmath,amssymb}
\usepackage[dvips]{graphicx, color}

\begin{document}
\title{\hskip 0.2cm\Large{$\mathbb{SYMMETRIES}$ \qquad}\qquad
\hskip 10.0cm\Large $\mathbb \&$\qquad\hskip 10.0cm 
{\Large{$\mathbb{THE\ RIEMANN\ HYPOTHESIS}$}}}
\author{\bf\large{Lin WENG}}
\date{(March 09, 2008)}
\maketitle
\noindent
{\bf Abstract} {\footnotesize Associated to classical semi-simple groups and their maximal parabolics are genuine zeta functions. Naturally related to Riemann's zeta and governed by symmetries, including that of Weyl, these zetas are expected to satisfy the Riemann hypothesis.}
\vskip 0.30cm
\noindent
For simplicity, $G$ here denotes a classical semi-simple algebraic group defined over the field $\mathbb Q$ of rationals.

With a fixed Borel, as usual, $\Delta_0$ stands for the corresponding collection of simple roots; $W$ the associatd Weyl group;  for a positive root $\alpha$, $\alpha^\vee$ the corresponding coroot; and  $\rho:=\frac{1}{2}\sum_{\alpha>0}\alpha$.
\vskip 0.30cm
\noindent
{\bf Definition 1.} The {\it period for $G$ over $\mathbb Q$} is defined by
$$\boxed{\omega_{\mathbb Q}^G(\lambda):=\sum_{w\in W}\Bigg(\frac{1}{\prod_{\alpha\in\Delta_0}\langle w\lambda-\rho,\alpha^\vee\rangle}\cdot\prod_{\alpha>0,w\alpha<0}\frac{\xi_{\mathbb Q}(\langle\lambda,\alpha^\vee\rangle)}
{\xi_{\mathbb Q}(\langle\lambda,\alpha^\vee\rangle+1)}\Bigg),\ \ \mathrm{Re}\,\lambda\in\mathcal C^+}$$ where $\mathcal C^+$ denotes the so-called positive chamber of $\frak a_0$, the space of 
characters associated to $(G,B)$, and 
$\xi_{\mathbb Q}(s)$ the completed Riemann zeta function.

For a fixed maximal parabolic subgroup $P$, it is well known that (the conjugation class of) $P$ corresponds to a simple root $\alpha_P\in\Delta_0$. Hence $\Delta_0\backslash\{\alpha_P\}=\{\beta_{1,P}, \beta_{2,P},\dots, \beta_{r-1,P}\},$ where $r=r(G)$ denotes the rank of $G$.
\vskip 0.30cm
\noindent
{\bf Definition 2.} The {\it period for $(G,P)$ over $\mathbb Q$} is defined by
$$\boxed{\begin{aligned}&\omega_{\mathbb Q}^{G/P}(\lambda_P):=\\
&\mathrm{Res}_{\langle \lambda-\rho,\beta_{r(G)-1,P}^\vee\rangle=0}\cdots
\mathrm{Res}_{\langle \lambda-\rho,\beta_{2,P}^\vee\rangle=0}\mathrm{Res}_{\langle \lambda-\rho,\beta_{1,P}^\vee\rangle=0}\,\Big(\omega_{\mathbb Q}^G(\lambda)\Big),\quad\mathrm{Re}\,\lambda_P\gg 0,\end{aligned}}$$
where with the constraint of taking residues along with  $(r-1)$ singular hyperplanes
$$\langle \lambda-\rho,\beta_{1,P}^\vee\rangle=0,\,
\langle \lambda-\rho,\beta_{2,P}^\vee\rangle=0,\,
\cdots,\,
\langle \lambda-\rho,\beta_{r(G)-1,P}^\vee\rangle=0,$$
there is only one variable, say a suitable $z_{i_0}$, left among $z_i$'s, 
re-scale it when necessary and rename it $\lambda_P$.
\vskip 0.30cm
Clearly, there is a minimal integer $I(G/P)$ and finitely many factors (depending on the choice of $\lambda_P$), 
$$\xi_{\mathbb Q}\Big(a_1^{G/P}\lambda_P+b_1^{G/P}\Big),\ \xi_{\mathbb Q}\Big(a_2^{G/P}\lambda_P+b_2^{G/P}\Big),\ \cdots,\ 
\xi_{\mathbb Q}\Big(a_{I(G/P)}^{G/P}\lambda_P+b_{I(G/P)}^{G/P}\Big),$$
such that {\it the product $\Big[\prod_{i=1}^{I(G/P)}\xi_{\mathbb Q}\Big(a_i^{G/P}\lambda_P+b_i^{G/P}\Big)\Big]\cdot \omega_{\mathbb Q}^{G/P}(\lambda_P)$ admits only finitely many singularities.}

\vskip 0.30cm
Similarly there is a minimal integer $J(G/P)$ and finitely many factors (depending on the choice of $\lambda_P$), 
$$\xi_{\mathbb Q}\Big(c_1^{G/P}\Big),\ \xi_{\mathbb Q}\Big(c_2^{G/P}\Big),\
\cdots,\
\xi_{\mathbb Q}\Big(c_{J(G/P)}^{G/P}\Big),$$ 
such that
{\it there are no factors of special $\xi_{\mathbb Q}$ values appearing at the denominators in  the product  $\Big[\prod_{i=1}^{J(G/P)}\xi_{\mathbb Q}\Big(c_i^{G/P}\Big)\Big]\cdot \omega_{\mathbb Q}^{G/P}(\lambda_P).$}
\vskip 0.30cm
\noindent
{\bf Definition 3.} (i) The {\it zeta function $\xi_{\mathbb Q;o}^{G/P}$ for $(G,P)$ over $\mathbb Q$} is defined by
$$\boxed{\begin{aligned}\xi_{\mathbb Q;o}^{G/P}\Big(s\Big)
:=&\Bigg[\prod_{i=1}^{I(G/P)}\xi_{\mathbb Q}\Big(a_i^{G/P}s+b_i^{G/P}\Big)\cdot\prod_{j=1}^{J(G/P)}\xi_{\mathbb Q}\Big(c_j^{G/P}\Big)\Bigg]\cdot \omega_{\mathbb Q}^{G/P}\Big(s\Big),\\
&\hskip 7.5cm\mathrm{Re}\,s\gg 0\end{aligned}}
$$

\noindent
{\bf Zeta Facts.} (1) {\it $\xi_{\mathbb Q;o}^{G/P}\Big(s\Big),\  \mathrm{Re}\,s\gg 0,$ is
a well-defined holomorphic function; admits a unique meromorphic continuation to the whole complex $s$-plane; and has only finitely many poles;} and

\noindent
(2) (Conjectural {\bf Functional Equation}) {\it There exists a constant $c_{G/P}\in\mathbb Q$ such that} $$\boxed{\xi_{\mathbb Q;o}^{G/P}\Big(-s+c_{G/P}\Big)=\xi_{\mathbb Q;o}^{G/P}\Big(s\Big).}$$

Obvious is (1). Rather complicated is (2), offering an additional symmetry.
\vskip 0.30cm
Classical symmetry $s\leftrightarrow 1-s$ for the standard functional equation  then leads to the following normalization.

\noindent
{\bf Definition 3.} (ii) The {\it zeta function $\xi_{\mathbb Q}^{G/P}\Big(s\Big)$ for $(G,P)$ over $\mathbb Q$} is defned by
$$\boxed{\xi_{\mathbb Q}^{G/P}\Big(s\Big):=\xi_{\mathbb Q;o}^{G/P}\Big(s+\frac{c_{G/P}-1}{2}\Big)}$$
\vskip 0.30cm
The most remarkable property shared by all these newly introduced zetas
is the following Zeta Fact about the uniformity of their zeros.
\vskip 0.20cm
\noindent
{\bf The Riemann Hypothesis\,$^{G/P}_{\mathbb Q}$.}

\noindent
$$\boxed{All\ zeros\ of\ the\ zeta\ function\ \xi_{\mathbb Q}^{G/P}\Big(s\Big)\ lie\ on\ the\ central\ line\
\mathrm{Re}\,s=\displaystyle{\frac{1}{2}}}$$

\eject
\centerline{\bf REFERENCES}
\vskip 0.30cm
\noindent
[Ar1] J. Arthur, A trace formula for reductive groups I. 
Terms associated to classes in $G({\mathbb Q})$. Duke Math. J. {\bf 45} (1978), 911--952
\vskip 0.20cm
\noindent
[Ar2] J. Arthur, A trace formula for reductive groups II: Applications of a 
truncation operator. Compositio Math. {\bf 40} (1980), no. 1, 87--121.
\vskip 0.20cm
\noindent
[Ar3] J. Arthur, A measure on the unipotent variety, Canad. J. Math {\bf 37}, (1985) 1237--1274
\vskip 0.20cm 
\noindent  
{[Co]} A. Connes, Trace formula in noncommutative geometry and the zeros 
of the Riemann zeta function.  Selecta Math. (N.S.) 5  (1999),  no. 1, 29--106.\vskip 0.20cm 
\noindent  
{[De]} C. Deninger, Motivic $L$-functions and regularized determinants,  Proc. Sympos. Pure Math, {\bf 55} Part 1, AMS (1994) 707-743
\vskip 0.20cm
\noindent
[D] B. Diehl, Die analytische Fortsetzung der Eisensteinreihe zur Siegelschen Modulgruppe, J. reine angew. Math., 317 (1980) 40-73
\vskip 0.20cm 
\noindent  
{[E]} H.M. Edwards, {\it Riemann's Zeta Function}, Dover (1974)
\vskip 0.20cm
\noindent 
[GS] G. van der Geer \& R. Schoof, Effectivity of Arakelov
Divisors and the Theta Divisor of a Number Field, Selecta Math. (N.S.)
{\bf 6} (2000), 377-398 
\vskip 0.20cm
\noindent
[Ha] T. Hayashi, Computation of Weng's rank 2 zeta function over an algebraic number field, J. Number Theory {\bf 125} (2007), no. 2, 473--527
\vskip 0.20cm 
\noindent
[Hu] J. Humphreys, {\it Introduction to Lie algebras and representations}, Springer-Verlag, 1972
\vskip 0.20cm 
\noindent  
{[Iw]}  K. Iwasawa, Letter to Dieudonn\'e, April 8, 1952,  Advanced
Studies in Pure Math. {\bf 21} (1992), 445-450
\vskip 0.20cm
\noindent 
[JLR] H. Jacquet, E. Lapid \& J. Rogawski,  Periods of automorphic forms. 
J. Amer. Math. Soc. 12 (1999), no. 1, 173--240
\vskip 0.20cm 
\noindent  
{[Ki]} H. Ki, All but finitely many non-trivial zeros of the approximations of 
the Epstein zeta function are simple and on the critical line. Proc. London Math. 
Soc. (3) 90 (2005), no. 2, 321--344. 
\vskip 0.20cm 
\noindent  
[KW] H.H. Kim \& L. Weng, Volume of truncated fundamental domains, Proc. AMS, Vol. 135 (2007) 1681-1688
\vskip 0.20cm
\noindent
[Laf] L. Lafforgue, {\it Chtoucas de Drinfeld et conjecture de 
Ramanujan-Petersson}. Asterisque No. 243 (1997)
\vskip 0.20cm
\noindent
[LS] J. Lagarias \& M. Suzuki, The Riemann Hypothesis for certain integrals of Eisenstein series, J. Number Theory, 118(2006), 98-122
\vskip 0.20cm
\noindent
[La] R. Langlands, {\it On the functional equations satisfied by 
Eisenstein series}, Springer LNM {\bf 544}, 1976
\vskip 0.20cm
\noindent
[La2] R. Langlands, Volume of fundamental domains for some 
arithmetical subgroups of Chevalley groups, Proc. Sympos. Pure Math. 9, AMS (1966) 143--148
\vskip 0.20cm 
\noindent  
[La3] R. Langlands,  {\it Euler products}, Yale Math. Monograph, Yale Univ.  1971
\vskip 0.20cm 
\noindent  
{[Mi]} H. Minkowski, {\it Geometrie der Zahlen}, Leipzig and Berlin, 1896
\vskip 0.20cm
\noindent
[MW] C. Moeglin \& J.-L. Waldspurger, {\it Spectral decomposition 
and Eisenstein series}. Cambridge Tracts in Math, {\bf 113}. 
Cambridge University Press, 1995
\vskip 0.20cm 
\noindent  
{[M]} D. Mumford, {\it Geometric Invariant Theory}, Springer-Verlag, 
 (1965)
 \vskip 0.20cm 
\noindent  
{[RR]} S. Ramanan \& A. Ramanathan,  Some remarks on the instability 
flag.  Tohoku Math. J. (2)  36  (1984),  no. 2, 269--291.
\vskip 0.20cm 
\noindent  
{[Ser]} J.-P. Serre, {\it Algebraic Groups and Class Fields}, GTM 117, 
Springer (1988)
\vskip 0.20cm
\noindent 
[Sie] C.L. Siegel, {\it Advanced analytic number theory},  T.I.F.R., Bombay, 1980
\vskip 0.20cm
\noindent   
{[St]} U. Stuhler,   Eine Bemerkung zur Reduktionstheorie quadratischer 
Formen, Arch. Math. (Basel) {\bf 27} (1976), no. 6, 604--610
\vskip 0.20cm
\noindent
[S] M. Suzuki, A proof of the Riemann Hypothesis for the Weng zeta function of rank 3 for the rationals, pp.175-200, in {\it Conference on L-Functions}, World Sci. 2007
\vskip 0.20cm
\noindent
[S2]  M. Suzuki, The Riemann hypothesis for Weng's zeta function of ${\rm Sp}(4)$ over $\mathbb Q$, preprint, 2008,
available at http://xxx.lanl.gov/abs/0802.0102
\vskip 0.32cm
\noindent
[SW] M. Suzuki \& L. Weng, private communications, Oct.-Dec., 2007
\vskip 0.20cm
\noindent
[SW2] M. Suzuki \& L. Weng, Zeta functions for $G_2$ and their zeros,  preprint, 2008, available at http://xxx.lanl.gov/abs/0802.0104
\vskip 0.20cm 
\noindent  
{[T]} J. Tate, Fourier analysis in number fields and Hecke's
zeta functions, Thesis, Princeton University, 1950
\vskip 0.20cm
\noindent
[W-1] L. Weng, Analytic truncation and Rankin-Selberg versus algebraic
truncation and non-abelian zeta, RIMS Kokyuroku, No.1324 (2003), 7-21
\vskip 0.20cm
\noindent
[W0] L. Weng, Non-abelian zeta function for function fields, Amer. J. Math 127 (2005), 973-1017
\vskip 0.20cm
\noindent
[W1] L. Weng, Geometric Arithmetic: A Program, in {\it Arithmetic Geometry and Number Theory}, pp. 211-390, World Sci. (2006)
\vskip 0.20cm
\noindent
[W2] L. Weng,  A Rank two zeta and its zeros, J of Ramanujan Math. Soc, 21 (2006), 205-266
\vskip 0.30cm
\noindent
[W3] L. Weng, A geometric approach to $L$-functions, in {\it Conference on L-Functions}, pp. 219-370, World Sci (2007)
\vskip 0.20cm
\noindent  
{[W4]} L. Weng, Zeta function for $Sp(2n)$, 
Appendix to [S2],  preprint, 2008, available at http://xxx.lanl.gov/abs/0802.0102
\vskip 0.20cm
\noindent
[Z] D. Zagier, The Rankin-Selberg method for auromorphic forms which are not of rapid decay, J. Fac. Sci., Univ. Tokyo, Sect. IA Math. 28 415-437 (1982)
\eject
\appendix
\section{Discovery of Zetas for $(G,P)/\mathbb Q$}
In this appendix, we expose some of the landmarks leading to the discovery of these elegant zetas introduced in the main text associated to reductive groups and their maximal parabolics over $\mathbb Q$. 

\bigskip
\noindent
{\bf Contents}
\vskip 0.50cm
\noindent
\begin{footnotesize}
{\bf A.1 High Rank Zeta Functions
\vskip 0.30cm
A.1.1 High rank zeta functions
\vskip 0.30cm
A.1.2 Relation with Eisenstein series
\vskip 0.30cm
A.1.3 $SL_2$, a toy model
\vskip 0.30cm
\noindent
A.2 General Periods
\vskip 0.30cm
A.2.1 Arthur's truncation and Eisenstein periods
\vskip 0.30cm
A.2.2 Rankin-Selberg \& Zagier method I: Sufficiently positive case
\vskip 0.30cm
A.2.3 Geo-arithmetical truncation and analytic truncation
\vskip 0.30cm
A.2.4 Rankin-Selberg \& Zagier method II: Semi-stable case
\vskip 0.30cm
A.2.5 Intertwining operator: Gindikin-Karpelevich formula
\vskip 0.30cm
A.2.6 Periods for $SL(n)$ over $\mathbb Q$: Weyl Symmetry
\vskip 0.30cm
\noindent
A.3 New Zetas for $SL(n)/\mathbb Q$
\vskip 0.30cm
A.3.1 Epstein, Koecher, Siegel zetas and Siegel-Eisenstein series
\vskip 0.30cm
A.3.2  Siegel's Eisenstein versus Langlands' Eisenstein
\vskip 0.30cm
A.3.3 New zetas: genuine but different
\vskip 0.30cm
A.3.4 Functional equation and the Riemann Hypothesis
\vskip 0.30cm
\noindent
A.4 Zetas for $(G,P)/\mathbb Q$
\vskip 0.30cm
A.4.1 From SL to Sp: analytic method adopted \& periods chosen
\vskip 0.30cm
A.4.2 $G_2$: maximal parabolics discovered
\vskip 0.30cm
A.4.3 Zetas for $(G,P)/\mathbb Q$: singular hyper-planes found
\vskip 0.30cm
\noindent
A.5 Conclusion remarks
\vskip 0.30cm
A.5.1 Analogue of high rank zetas
\vskip 0.30cm
A.5.2 $T$-version
\vskip 0.30cm
A.5.3 Where lead to}
\end{footnotesize}
\eject
\subsection{High Rank Zeta Functions}

\noindent
{\bf A.1.1 High Rank Zeta Functions}

\smallskip
\noindent
Let $F$ be a number field with $\mathcal O_F$ the ring of integers. Denote by $\Delta_F$ the discriminant of $F$. Fix a positive integer $r$. Then by definition, an $\mathcal O_F$-lattice $\Lambda$ of rank $r$ is a pair $(P,{\bf \rho})$ consisting of an $\mathcal O_F$-projective module $P$ of rank $r$ and a metric
${\bf \rho}:=(\rho_{\sigma:\mathbb R},\rho_{\tau:\mathbb C})$ on $\Big(\mathbb R^{r_1}\times\mathbb C^{r_2}\Big)^r=\Big(\mathbb R^r\Big)^{r_1}\times\Big(\mathbb C^r\Big)^{r_2}$. Here as usual, we denote by $r_1$ and $r_2$ the number of real embeddings $\sigma:F\hookrightarrow\mathbb R$  and the numbers of complex embeddings  $\tau:F\hookrightarrow \mathbb C$, respectively. (Recall that by a standard result, $P$ is isomorphic to $\mathcal O_F^{\oplus(r-1)}\oplus\frak a$ for a suitable fractional ideal $\frak a$ of $F$. Thus via the natural inclusion $\mathcal O_F^{\oplus(r-1)}\oplus\frak a\hookrightarrow F^r\hookrightarrow \Big(\mathbb R^{r_1}\times\mathbb C^{r_2}\Big)^r$, we may view $P$ as a discrete subgroup of the matrized space
$\Big(\mathbb R^{r_1}\times\mathbb C^{r_2}\Big)^r$. Fix it.)
It is well known that the quotient space $\Big(\mathbb R^{r_1}\times\mathbb C^{r_2}\Big)^r/\Lambda$ is compact. Call its volume the (co-)volume of $\Lambda$ and denote it by $\mathrm{Vol}(\Lambda)$.
By definition, a lattice $\Lambda$ is called {\it semi-stable} if for all $\mathcal O_F$-sublattices $\Lambda_1$, we have
$$\mathrm{Vol}(\Lambda_1)^{\mathrm{rk}(\Lambda)}\geq \mathrm{Vol}(\Lambda)^{\mathrm{rk}(\Lambda_1)}.$$
Denote by $\mathcal M_{F,r}$ the moduli space of semi-stable $\mathcal O_F$-lattices of rank $r$. (For details, see e.g., [W1-3],  [Gr1,2], [St1,2].) This is the first ingredient needed to introduce  high rank zetas for $F$.  In particular, we know the following
\vskip 0.30cm
\noindent
{\bf Fact A}. ([W1-3]) (1) {\it There is a natural decomposition $\mathcal M_{F,r}=\cup_{T\in\mathbb R_{>0}}\mathcal M_{F,r}[T]$ where, $\mathcal M_{F,r}[T]$ denotes the moduli space of semi-stable $\mathcal O_F$-lattices of rank $r$ and of volume $T$;}

\noindent
(2) {\it $\mathcal M_{F,r}[T]$ is compact;} and

\noindent
(3) {\it There are natural measures $d\mu$ and $d\mu_0$ on $\mathcal M_{F,r}$ and on $\mathcal M_{F,r}[T]$ respectively such that with respect to the decomposition (1), we have $d\mu=d\mu_0\times\frac{dT}{T}$.}
\vskip 0.30cm
The second ingredient needed is a good geo-arithmetical cohomology. For this, we define the 0-th cohomology group $H^0(F,\Lambda)$ of an $\mathcal O_F$-lattice $\Lambda$ to be the the lattice $\Lambda$ itself,
and the 1-st cohomology group $H^1(F,\Lambda)$ to be the compact quotient group $\Big(\mathbb R^{r_1}\times\mathbb C^{r_2}\Big)^r/\Lambda$. Consequently, we have the following
Pontryagin duality for them:

\noindent
{\bf Topological Duality}. $\widehat{H^1(F,\Lambda)}\simeq H^0(F,\omega_F\otimes\Lambda^\vee).$

\noindent
Here for a locally compact group $G$, denote by $\widehat G$
its Pontryagin dual, $\omega_F$ denotes the differential lattice, i.e., the lattice whose module part is simply the module of differentials of $F$, while whose metric is the standard one on $\mathbb R^{r_1}\times\mathbb C^{r_2}$. As such, following Tate ([T] and [W]), we then can use Fourier analysis to count our $H^0(F,\Lambda)$ and $H^1(F,\Lambda)$. For example, each element ${\bf x}\in H^0(F,\Lambda)$ is counted with the weight of Gaussian distribution
$\sum_{\sigma:\mathbb R}e^{-\pi\|{\bf x}\|_{\rho_\sigma}^2}+
\sum_{\tau:\mathbb C}e^{-2\pi\|{\bf x}\|_{\rho_\tau}^2}$ and accordingly define $h^0(F,\Lambda)$ to be the logarithm of this count. (See also [GS].) 
Particularly, with such $h^0$ and $h^1$ for a lattice $\Lambda$, by using the above topological duality and the Poisson summation formula,
 then we obtain the following
\vskip 0.30cm
\noindent
{\bf Fact B.} ([W1-3]) {\it Let $\Lambda$ be an $\mathcal O_F$-lattice of rank $r$. Then} 

\noindent
(1) ({\bf Duality}) $h^1(F,\Lambda)=h^0(F,\omega_F\otimes\Lambda^\vee);$ and

\noindent
(2) ({\bf Riemann-Roch Theorem}) $h^0(F,\Lambda)-h^1(F,\Lambda)
=\mathrm{deg}\,(\Lambda)-\frac{r}{2}\log|\Delta_F|.$ {\it Here $\mathrm{deg}\,\Lambda$ denotes the Arakelov degree of $\Lambda$.}

\noindent
(For the reader who does not know Arakelov degree, recall then the following weak result\newline {\bf Arakelov-Riemann-Roch Theorem}:
$-\log\,\mathrm{Vol}\,(\Lambda)=\mathrm{deg}\,(\Lambda)-\frac{r}{2}\log|\Delta_F|.$)
\vskip 0.30cm
With all this, then we are ready to introduce the following
\vskip 0.30cm
\noindent
{\bf Definition.} ([W1,3]) For an algebraic number field $F$ and a positive integer $r$, define 
its {\it rank $r$ zeta function} by
$$\xi_{F,r}(s):=\Big(|\Delta_F|\Big)^{\frac{r}{2}s}\int_{\mathcal M_{F,r}}\Big(e^{h^0(F,\Lambda)}-1\Big)\cdot\Big(e^{-s}\Big)^{\mathrm{deg}(\Lambda)}\,d\mu(\Lambda),\qquad\mathrm{Re}\,(s)>1.$$

From the definition, by Fact A for moduli spaces and Fact B on Duality and the Riemann-Roch for geo-arithmetic cohomologies, totologically, we have the following
\vskip 0.30cm
\noindent
{\bf Fact C.} ([W1,3]) (0) (Iwasawa) {\it $\xi_{F,1}(s)\buildrel{\cdot}\over{=}\xi_F(s)$,  the completed Dedekind zeta function;}

\noindent
(1) ({\bf Mero Extension}) {\it $\xi_{F,r}(s)$ is well-defined and admits a meromorphic continuation to the whole complex $s$-plane;}

\noindent
(2) ({\bf Functional Equation}) $\xi_{F,r}(1-s)=\xi_{F,r}(s)$; and

\noindent
(3) ({\bf Singularities}) {\it There are only two singularities,  i.e.,   simple poles at $s=0,\,1$  with the residue $\mathrm{Res}_{s=1}\xi_{F,r}(s)=\mathrm{Vol}\,\mathcal M_{F,r}\Big(|\Delta_F|^{\frac{r}{2}}\Big).$}
\vskip 0.30cm
\noindent
{\bf A.1.2 Relation with Eisenstein Periods}

\smallskip
\noindent
We next give a relation between our high rank zetas and what we call  Eisenstein periods.
The point here is instead of working over $\mathcal M_{F,r}$, we fix a volume so as to work over the {\it compact} subspace $\mathcal M_{F,r}[|\Delta_F|^{r/2}]$ and hence deduce the desired relation via Mellin transform. This goes as follows.

From now on, for simplicity, we work over the field $\mathbb Q$ of rationals. Accordingly, the rank $r$ zeta function 
$\xi_{{\mathbb  Q},r}(s)$ of ${\mathbb  Q}$ is given by
$$\xi_{{\mathbb  Q},r}(s)=
\int_{{\mathcal  M}_{{\mathbb  Q},r}}\left(e^{h^0({\mathbb  Q},\Lambda)}-1
\right)
\cdot \big(e^{-s}\big)^{\mathrm{deg}(\Lambda)}\, d\mu(\Lambda),\qquad 
\mathrm{Re}(s)>1,$$
where $h^0({\mathbb  Q},\Lambda):=\log\left(\sum_{x\in \Lambda}
\exp\big(-\pi|x|^2\big)\right)$ and $\mathrm{deg}(\Lambda)=-\log\,
\mathrm{Vol}\big(\mathbb R^r/\Lambda\big).$

Decompose according to their volumes, 
${\mathcal  M}_{{\mathbb  Q},r}=\cup_{T>0}{\mathcal  M}_{{\mathbb  Q},r}[T]$,
and there is a natural morphism
${\mathcal  M}_{{\mathbb  Q},r}[T]\to {\mathcal  M}_{{\mathbb  Q},r}[1],
\ \Lambda\mapsto T^{1\over r}\cdot\Lambda.$
Consequently, $$\begin{aligned} \xi_{{\mathbb  Q},r}(s)=&
\int_{\cup_{T>0}{\mathcal  M}_{{\mathbb  Q},r}[T]}\left( e^{h^0({\mathbb  Q},
\Lambda)}-1
\right)\cdot \big(e^{-s}\big)^{\mathrm{deg}(\Lambda)}\, d\mu(\Lambda)\\
=&\int_0^\infty T^s{{dT}\over T}
\int_{{\mathcal  M}_{{\mathbb  Q},r}[1]}\left( 
e^{h^0({\mathbb  Q},T^{1\over r}\cdot \Lambda)}-1\right)
\,d\mu(\Lambda).\end{aligned}$$
But $h^0({\mathbb  Q},T^{1\over r}\cdot \Lambda)=
\log\left(\sum_{x\in \Lambda}
\exp\big(-\pi|x|^2\cdot T^{2\over r}\big)\right)$. By applying the Mellin transform, we have
$$\xi_{{\mathbb  Q},r}(s)={r\over 2}\cdot\pi^{-{r\over 2}\, s}
\Gamma\Big({r\over 2}\, s\Big)\cdot\int_{{\mathcal  M}_{{\mathbb  Q},r}[1]}
\Big(\sum_{x\in\Lambda\backslash\{0\}}|x|^{-rs}\Big)\, d\mu_0(\Lambda).$$
Accordingly, introduce the completed Epstein zeta function for
$\Lambda$ by 
$$\hat E(\Lambda;s):=\pi^{-s}\Gamma(s)\cdot \sum_{x\in \Lambda\backslash \{0\}}
|x|^{-2s}.$$ We then arrive at

\noindent
{\bf Fact D.} ([W1-3]) (Eisenstein series and high rank zetas)
{\it $$\xi_{\mathbb Q,r}(s)=\frac{r}{2}\int_{\mathcal M_{\mathbb Q,r}[1]}\widehat E(\Lambda,\frac{r}{2}s)d\mu_0(\Lambda),\qquad\mathrm{Re}(s)>1.$$}

\vskip 0.30cm
\noindent
{\bf A.1.3 $SL(2)$: A Toy Model}

\smallskip
\noindent
To indicate basic ideas clearly, we first give some details on the rank two zeta $\xi_{\mathbb Q,2}(s)$.

Consider the action of $\mathrm{SL}(2,{\mathbb  Z})$ on the upper half plane 
${\mathcal  H}(=SL(2,\mathbb R)/SO(2))$.
Then we obtain a standard  \lq fundamental domain\rq$\ $  
$D=\{z=x+iy\in {\mathcal  H}:|x|\leq {1\over 2},y>0,x^2+y^2\geq 1\}.$ 
Recall also the completed standard Eisenstein series 
$$\hat E(z;s):=\pi^{-s}\Gamma(s)\cdot\sum_{(m,n)\in {\mathbb  Z}^2
\backslash \{(0,0)\}}{{y^{s}}\over {|mz+n|^{2s}}}.$$
Naturally, we are led to considering the integral
$\int_D\hat E(z,s){{dx\,dy}\over {y^2}}.$
However, this integration diverges. Indeed,
near the only cusp $y=\infty$, by the Chowla-Selberg formula,
$\hat E(z,s)$ has the  Fourier expansion
$$\hat E(z;s)=\sum_{n=-\infty}^\infty a_n(y,s)e^{2\pi i nx}$$ with 
$$a_n(y,s)=\begin{cases} \xi(2s)y^s+\xi(2-2s)y^{1-s},&  \text{if $n=0$;}\\
2|n|^{s-{1\over 2}}\sigma_{1-2s}(|n|){\sqrt y}K_{s-{1\over 2}}(2\pi|n|y),& 
\text{if $n\not=0$,} \end{cases}$$
where $\xi(s)$ is the completed Riemann zeta function,
$\sigma_s(n):=\sum_{d|n}d^s$, and 
$K_s(y):={1\over 2}\int_0^\infty e^{-y(t+{1\over t})/2}t^s{{dt}\over t}$ 
is the K-Bessel function. Moreover, 
$$|K_s(y)|\leq e^{-y/2}K_{\mathrm{Re}(s)}(2),\ \mathrm{if}\  y>4,\qquad
\mathrm{and}\qquad K_s=K_{-s}.$$ 
So $a_{n\not=0}(y,s)$ decay exponentially, and the constant term $a_0(y,s)$, being of slow growth, is  problematic.

Therefore, to introduce a meaningful integration from the original ill-defined one, 
we need cut off the slow growth part. There are two ways to do so: 
one is geometrical and hence rather direct and simple; the other is 
analytical, and hence rather technical and traditional, dated back to 
Rankin-Selberg.
\vskip 0.30cm
\noindent
(a) {\bf Geometric Truncation}

Draw a horizontal line $y=T\geq 1$ and set
 $$D_T=\{z=x+iy\in D:y\leq T\},\qquad D^T=\{z=x+iy\in D:y\geq T\}.$$ 
Then $D=D_T\cup D^T$.
Introduce a well-defined integration
$$I^{\mathrm{Geo}}_T(s):=\int_{D_T}\hat E(z,s)\,{{dx\,dy}\over {y^2}}.$$

\noindent
(b) {\bf Analytic Truncation}

Define a truncated Eisenstein series $\hat E_T(z;s)$ by
$$\hat E_T(z;s):=\begin{cases} \hat E(z;s),&\text{if $y\leq T$;}\\
\hat E(z,s)-a_0(y;s),&\text{if $y>T$.}\end{cases}$$
Introduce a well-defined integration
$$I_T^{\mathrm{Ana}}(s):=\int_D\hat E_T(z;s)\,{{dx\,dy}\over {y^2}}.$$

With this, from the Rankin-Selberg method,  we have the following:
\vskip 0.30cm
\noindent
{\bf Fact E.} (See e.g., [Z]) (Analytic Truncation=Geometric Truncation in Rank 2) 
{\it $$I_T^{\mathrm{Geo}}(s)=\frac{\xi(2s)}{s-1}\cdot T^{s-1}-\frac{\xi(2s-1)}{s}\cdot
T^{-s}=I_T^{\mathrm{Ana}}(s).$$}

Each of the above two integrations has its own merit: for the geometric one, 
we keep the Eisenstein series unchanged, while for the analytic one, we keep
the original fundamental domain of ${\mathcal  H}$ under 
$\mathrm{SL}(2,{\mathbb  Z})$
as it is.

Note that a particular nice point about the fundamental domain is that it admits a
modular interpretation. Thus it would be very nice if we could on the one hand 
keep the Eisenstein series unchanged, while on the other hand offer some integration
domains which appear naturally in certain moduli problems. This is essential the idea of introducing $\mathcal M_{F,r}\Big(|\Delta_F|^{\frac{r}{2}}\Big)$, the first key ingredient for high rank zetas.
\vskip 0.30cm
\noindent
(c) {\bf  Arithmetic Truncation}

Now we explain why the above discussion and Rankin-Selberg method have 
anything to do with our high rank zeta functions. For this, we introduce 
yet another truncation, the geo-arithmetic one using stability.

So back to the moduli space  of rank 2 lattices of volume 1 over 
${\mathbb  Q}$.
Then classical reduction theory gives a  natural map from this moduli
space to the fundamental domain $D$ above: 
For any lattice $\Lambda$ in $\mathbb R^2$, fix ${\bold x}_1\in \Lambda$ such that its length  
gives the first Minkowski minimum $\lambda_1$ of $\Lambda$.
Then via rotation, we may assume that ${\bold x}_1=(\lambda_1,0)$. 
Further, from the reduction theory 
${1\over {\lambda_1}}\Lambda$ may be viewed as the lattice of the volume
$\lambda_1^{-2}=y_0$ generated by $(1,0)$ and $\omega=x_0+iy_0\in D$. 
That is to say, 
the points in $D_T$ constructed in (a) above are in one-to-one corresponding to
rank two lattices of volume one whose
first Minkowski minimum $\lambda$, satisfying $\lambda_1^{-2}\leq T$, i.e,
$\lambda_1\geq T^{-{1\over 2}}$. Set 
${\mathcal  M}_{{\mathbb  Q},2}^{\leq {1\over 2}\log T}[1]$  be the moduli 
space of
rank 2 lattices $\Lambda$ of volume 1 over ${\mathbb  Q}$ all of whose sublattices 
$\Lambda_1$ of rank 1 have degrees 
$\leq {1\over 2}\log T$. With this discussion, we have the following
\vskip 0.30cm
\noindent
{\bf Fact F.} ([W1-3]) (Geometric Truncation = Arithmetic Truncation) \newline{\it 
There is a natural (quasi) one-to-one, onto morphism 
$${\mathcal  M}_{{\mathbb  Q},2}^{\leq {1\over 2}\log T}[1]\simeq D_T.$$
In particular, $${\mathcal  M}_{{\mathbb  Q},2}^{\leq 0}[1]=
{\mathcal  M}_{{\mathbb  Q},2}[1]\simeq D_1.$$}

Consequently,  we have the following
\vskip 0.30cm
\noindent
{\bf Example in Rank 2}. ([W1-3]) {\it  
$\displaystyle{\xi_{{\mathbb  Q},2}(s)={{\xi(2s)}\over {s-1}}-{{\xi(2s-1)}\over {s}}.}$}

\subsection{Periods}

\noindent
{\bf A.2.1 Arthur's Truncation and Eisenstein Periods}

\smallskip
\noindent
Recall that the upper half plane $\mathcal H$ admits the following group theoretic interpretation
$SL(2,\mathbb Z)\backslash SL(2,\mathbb R)/SO(2)$. Thus for high rank zeta functions, we then naturally shift to $G=SL(n)$, or more generally, any split group $G$.

Fix  a parabolic subgroup $P$ of $G$ with Levi decomposition $P=MN$, denote by $\frak a_P$ the complexification of the space of characters associated to $P$. In particular, denote by $\frak a_0$ the one for the Borel.
Denote by $\Delta_0$ the associated collection of simple roots. By definition, an element $T\in\frak a_0$ is said to be {\it sufficiently regular}, or sufficiently positive, and denoted by $T\gg0$ if for all $\alpha\in\Delta_0$ $\langle \alpha,T\rangle\gg 0$ are large enough. Fix such a $T$. 

Let $\phi:G(\mathbb Z)\backslash G(\mathbb R)/K\to \mathbb C$ be a smooth function where $K$ is a maximal compact subgroup of $G(\mathbb R)$. We define 
{\it Arthur's analytic truncation} $\wedge^T\phi$ (for $\phi$ with respect to the parameter $T$) to be the function on $G(\mathbb Z)\backslash G(\mathbb R)/K$ given by $$\Big(\wedge^T\phi\Big)(Z)
:=\sum_{P:\mathrm{standard}}(-1)^{\mathrm{rank}(P)}\sum_{\delta\in P(\mathbb Z)\backslash G(\mathbb Z)}\phi_P(\delta g)\cdot\hat\tau_P\big(H_P(\delta g)-T\big),$$ where $\phi_P:=\int_{N(\mathbb R)/N(\mathbb R)\cap SL(n,\mathbb Z)}f(xn)\,dn$ denotes the constant term of $\phi$ along with the standard parabolic subgroup $P$, $\hat\tau_P$ is the characteristic function of the so-called positive cone in $\frak a_P$, and $H_P(g):=\log_Mm_P(g)$ is an elelemnt in $\frak a_P$. (For unknown notation, all standard, see e.g., [Ar1,2], [JLR], or  [W-1,3].)

Fundamental properties of Arthur's truncation may be summarized as:
\vskip 0.30cm
\noindent
{\bf Fact G.} ([Ar1,2] \&/or [OW]) {\it For  a sufficiently positive $T$ in $\frak a_0$, 
 we have
\begin{description}
\item [(1)] $\wedge^T\phi$ is rapidly decreasing, if $\phi$ is an automorphic 
form on $G(\mathbb Z)\backslash G(\mathbb R)/K$;
\item [(2)]  $\wedge^T\circ \wedge^T=\wedge^T$;
\item [(3)] $\wedge^T$ is self-adjoint; and
\item [(4)]  ([Ar4]) $\wedge^T{\bold 1}$ is a characteristic function of a compact 
subset of $G(\mathbb Z)\backslash G(\mathbb R)/K$.
\end{description}}

Denote by $\frak F(T)$ the compact subset of $G(\mathbb Z)
\backslash G(\mathbb R)/K$ whose characteristic function is given by 
$\wedge^T\bold 1$ by (4).
\vskip 0.30cm
\noindent
{\bf Corollary.} ([W1,3]) {\it Let $T\gg0$ be a fixed element in $\frak a_0$. 
If $\phi$ is an automorphic form on $G(\mathbb Z)\backslash G(\mathbb R)/K$,
$$\int_{G(\mathbb Z)\backslash G(\mathbb R)/K}\wedge^T\phi(g)\,
dg=\int_{\frak F(T)}\phi(g)\,dg.
$$}

We call the above integration the {\it Arthur periods} associated to $\phi$. In most of applications, the following special class, called Eisenstein periods, plays a key role.

Recall that if $\varphi$ is an $M$-level automorphic form, then we may form the associated Eisenstein series $E^{G/P}(\varphi,\lambda)(g)=E(\varphi,\lambda)(g)$ as follows:
$$E(\varphi,\lambda)(g):=\sum_{P(\mathbb Z)\backslash G(\mathbb Z)}
m_P(\delta g)^{\lambda+\rho_P}\cdot\phi(\delta g),\quad\mathrm{Re}\,\lambda\in\mathcal C_P^+$$ where $\mathcal C_P^+$ denotes the positive chamber in $\frak a_P$. By definition, the {\it Eisenstein period}
is the integration $$\int_{G(\mathbb Z)\backslash G(\mathbb R)/K}\wedge^TE(\varphi,\lambda)(g)\,
dg=\int_{\frak F(T)}E(\varphi,\lambda)(g)\,dg.$$
\noindent
(Here we use a normalization
for the Eisenstein series as usual, i.e., shifting the variable from $\lambda$ to $\lambda+\rho_P$, 
so that the convergence region is simply the positive chamber.)

\vskip 0.30cm
\noindent
{\bf A.2.2 Rankin-Selberg \& Zagier Method \newline
\hskip 1.15cm I: Sufficiently Positive Case}

\smallskip
\noindent
In general, it is {\it very difficult}, in fact, {\it quite impossible,}  to calculate Eisenstein period precisely. However, if the original automorphic function (in defining the Eisenstein series used) is {\it cuspidal}, this can be evaluated. 
This is a result due to Jaquet-Lapid-Rogowski (see e.g., [JLR]), which itself may be viewed as an advanced version of the so-called Rankin-Selberg \& Zagier method.  (See also section 4.2 [W0] for our own solution, which was quite similar and was independently written before we knew [JLR].) 
In particular, for constant function $\bold 1$ over the Borel, and the associated Eisenstein series by $E({\bf 1};\lambda;g)$,
we have the following:
\vskip 0.30cm
\noindent
{\bf Fact E}$'$. ([JLR], [W0]) {\it Assume that $T$ is sufficiently positive, then the Eisensetin period $\int_{G(\mathbb Z)\backslash G(\mathbb R)} \wedge^T E({\bf 1};\lambda;g)\, dg$ is given by
$$\int_{G(\mathbb Z)\backslash G(\mathbb R)} \wedge^T E({\bf 1};\lambda;g)\, dg=v\sum_{w\in W} \frac {e^{\langle w\lambda-\rho, T\rangle}}{\prod_{\alpha\in \Delta_0} \langle w\lambda-\rho, \alpha^{\vee}\rangle}\cdot  M(w,\lambda)$$ 
where $v=\mathrm{Vol}\big(\{\sum_{\alpha\in \Delta_0} a_{\alpha}
\alpha^{\vee}: 0\leq a_{\alpha}<1\}\big)$, $W$ denotes the Weyl group, 
$\Delta_0$ the set of simple roots, $\alpha^\vee$ the co-root associated to $\alpha$, and $M(w,\lambda)$ denotes the assosciated intertwining operator.}
\vskip 0.30cm
\noindent
{\bf A.2.3 Geo-Arithmetic Truncation and Analytic Truncation}

\smallskip
\noindent
In algebraic geometry, or better in Geometric Invariant Theory, 
a fundamental principle, which we call {\it the Micro-Global Principle}, 
claims that if a point is not GIT stable then there exists
a parabolic subgroup which destroys the corresponding stability.

Here even we do not have a proper definition of GIT stability for lattices, 
in terms of intersection stability, an analogue of the Micro-Global Principle 
holds. To see this, we go as follows 
(and for our own convenience, we adopt an adelic language when necessary).

For $g=g(\Lambda)\in G(\mathbb A)$, denote its associated 
lattice by $\Lambda^g$, and its induced filtration from $P$ by
$$0=\Lambda_0^{g,P}\subset \Lambda_1^{g,P}\subset\cdots\subset 
\Lambda_{|P|}^{g,P}=\Lambda^g.$$ 
(Recall that all lattices can be obtained in this manner, and that for a fixed lattice, its associated fiber in $G(\mathbb A)$ is compact.) 
Assume that $P$ corresponds to the partition $I=(d_1,d_2,\ldots,d_{n=:|P|})$. 
Consequently, we have $$\mathrm{rk}(\Lambda_i)=r_i:=d_1+d_2+\cdots+d_i,\qquad 
\mathrm{for}\ i=1,2,\ldots,|P|.$$ 
Define the polygon 
$p_P^g=p_P^{\Lambda^g}:[0,r]\to\mathbb R$ of $\Lambda=\Lambda^g$ with respect to $P$ by

\noindent
(1) $p_P^g(0)=p_P^g(r)=0$;

\noindent
(2) $p_P^g$ is affine on $[r_i,r_{i+1}]$, $i=1,2,\ldots,|P|-1$; and

\noindent
(3) $p_P^g(r_i)=\mathrm{deg}\big(\Lambda_i^{g,P}\big)-r_i\cdot
\frac{\mathrm{deg}\big(\Lambda^{g}\big)}{r},\ i=1,2,\ldots,|P|-1$.

\noindent
Note that if the volume of $\Lambda$ is assumed to be one, then 
(3) is equivalent to 

\noindent
(3)$'$ $p_P^g(r_i)=\mathrm{deg}\big(\Lambda_i^{g,P}\big),\ i=1,2,\ldots,|P|-1$.

Based on stability, we may introduce a more general 
geometric truncation for the space of lattices. For this we start with the following easy statement:

\noindent
{\it For a fixed ${\mathcal O}_F$-lattice 
$\Lambda$,}
$\Big\{\mathrm{Vol}(\Lambda_1):\Lambda_1\subset\Lambda\Big\}\subset 
{\mathbb R}_{\geq 0}$ {\it is  discrete and bounded from
below.}

As a direct consequence, we have the following

\noindent
{\bf Fact H.} ([W1,3]) ({Canonical Filtration}) {\it For an ${\mathcal O}_F$-lattice $\Lambda$,
there exists a unique 
filtration, called the canonical filtration of $\Lambda$, of proper sublattices
$$0=\Lambda_0\subset\Lambda_1\subset\cdots \subset\Lambda_s=\Lambda$$ 
such that}

\noindent
(1) {\it for all $i=1,\cdots, s$, $\Lambda_i/\Lambda_{i-1}$ 
is semi-stable}; and

\noindent
(2) {\it for all} $j=1,\cdots s-1$,
$$\Big(\mathrm{Vol}(\Lambda_{j+1}/\Lambda_j)\Big)^{
\mathrm{rk}(\Lambda_j/\Lambda_{j-1})}>
\Big(\mathrm{Vol}(\Lambda_{j}/\Lambda_{j-1})\Big)^{
\mathrm{rk}(\Lambda_{j+1}/\Lambda_{j})}.$$

Accordingly, for an $\mathcal O_F$-lattice $\Lambda$ 
with the associated canonical filtration, (an analogue of the Harder -Narasimhan-Langton filtration for vector bundles over Riemann surfaces [HN],)
$$0=\overline\Lambda_0\subset\overline\Lambda_1\subset\cdots 
\subset\overline\Lambda_s=\Lambda$$
define the associated {\it canonical polygon} 
$\overline p_\Lambda:[0,r]\to {\mathbb R}$ by the following conditions:
\begin{description}
\item[(1)] $\overline p_\Lambda(0)=\overline p_\Lambda(r)=0$;

\item[(2)] $\overline p_\Lambda$ is affine over the closed interval
$[\mathrm{rk}\overline\Lambda_i,\mathrm{rk}\overline\Lambda_{i+1}]$; and

\item[(3)] $\overline p_\Lambda(\mathrm{rk}\overline\Lambda_i)
=\mathrm{deg}(\overline\Lambda_i)-\mathrm{rk}(\overline\Lambda_i)
\cdot\frac{\mathrm{deg}(\overline\Lambda)}{r}.$
\end{description} 

Let now $p,q:[0,r]\to {\mathbb R}$ be two polygons 
such that $p(0)=q(0)=p(r)=q(r)=0$. Then, 
we say $q$ is {\it bigger than} $p$ {\it with respect to} $P$ and denote it by 
$q>_Pp$, if $q(r_i)-p(r_i)>0$ for all $i=1,\ldots,|P|-1.$ (See e.g., [Laf].)
Introduce also  the characteristic function $\bold 1(\overline p^*\leq p)$ by
$$\bold 1(\overline p^g\leq p)
=\begin{cases} 1,&\text{if $\overline p^g\leq p$;}\\
0,&\text{otherwise}.\end{cases}$$ Here $\overline p^g$ denotes 
the canonical polygon for the lattice corresponding to $g$. 
Recall that for a parabolic subgroup 
$P$, $p_P^g$ denotes the polygon induced by $P$ for (the lattice corresponding 
to) the element $g\in G(\mathbb A)$.
\vskip 0.30cm
\noindent
{\bf Fact I.} ([W1,3]) (Fundamental Relation) 
{\it For a fixed convex polygon $p:[0,r]\to 
{\mathbb R}$ such that $p(0)=p(r)=0$, we have
$$\bold 1(\overline p^g\leq p)=\sum_{P:\, \mathrm{standard\, 
parabolic}}(-1)^{|P|-1}\sum_{\delta\in
P(F)\backslash G(F)}\bold 1(p_P^{\delta g}>_Pp).$$}

\noindent
{\it Remarks.} (1) This is an arithmetic analogue of a result of Lafforgue ([Laf]) for vector bundles over function fields. 

\noindent
(2) The right hand side may be naturally decomposite into 
two parts according to whether $P=G$ or not.
In such  away,  the right hand side becomes
$$\bold 1_G-\sum_{P:\, \mathrm{proper\,standard\, parabolic}}(-1)^{|P|-1}
\cdots.$$ This then exposes two aspects of our geometric truncation:
First of all, if a lattice is not stable, then there are parabolic
subgroups which take the responsibility; Secondly, each parabolic subgroup
has its fix role -- Essentially, they should be counted only once 
each time. In other words, if more are substracted, then we need  
add one fewer back to make sure the whole process is not overdone.
\vskip 0.30cm
From (2) above, it is clear that the geo-arithmetical truncation defined using $\bold 1(\overline p^g\leq p)$, or simply using stability, has the same strucrure as that for analytic truncations. Next, we want to give a precise relation between these two truncations, so that analytic methods created in the study of trace formula
can be employed in the study of our high rank zetas.
\vskip 0.30cm
Recall that a polygon $p:[0,r]\to\mathbb R$ is called {\it normalized} 
if $p(0)=p(r)=0$. For a (normalized) polygon $p:[0,r]\to\mathbb R$,
define the associated (real) character $T=T(p)\in\frak a_0$ of 
$M_0$ (the Levi for the Borel) by the condition that $$\alpha_i(T)=\Big[p(i)-p(i-1)\Big]-
\Big[p(i+1)-p(i)\Big]$$ 
for all $i=1,2,\ldots, r-1,$ where $\alpha_i=e_i-e_{i+1}\in\Delta_0$ denote simple roots. 
As such, one checks that
$$T(p)=\Big(p(1)-p(0), \cdots,
\cdots, p(i)-p(i-1),\dots,p(r)-p(r-1)\Big).$$

Set also ${\bf 1}(p_P^*>_Pp)$ to be the characteristic function 
of the subset
of $g$'s such that $p_P^g>_Pp$. Then we have the following
\vskip 0.30cm
\noindent
{\bf Fact J.} ([W1,3]) (Micro Bridge) {\it For a fixed convex normalized polygon 
$p:[0,r]\to\mathbb R$, and $g\in SL_r(\mathbb A)$, with respect to
 any parabolic subgroup $P$, we have
 $$\hat\tau_P\Big(-H_0(g)-T(p)\Big)=\bold 1\Big(p_P^g>_Pp\Big).$$}

\vskip 0.30cm
With this micro bridge, we are ready to expose a beautiful 
intrinsic relation between our geo-arithmetic truncation using stability and analytic truncations.
\vskip 0.30cm
\noindent
{\bf Fact J.} ([W1,3]) (Global Bridge) {\it For a fixed normalized convex polygon $p:[0,r]\to 
{\mathbb R}$, let $$T(p):=\Big(p(1), p(2)-p(1),
\dots, p(i)-p(i-1),\dots,p(r-1)-p(r-2),-p(r-1)\Big)$$ be the
associated vector in $\frak a_0$. If $T(p)$ is sufficiently positive,
then $$\bold 1(\overline p^g\leq p)=\Big(\Lambda^{T(p)}\bold 1\Big)(g).$$} 

In particular, by Facts G, I, and J, we arrive at the following analytic interpretation of the moduli space of semi-stable lattices.
\vskip 0.30cm
\noindent
{\bf Fact G-I-J.} ([W1,3]) $\frak F(0)=\mathcal M_{\mathbb Q,r}[1].$
\vskip 0.30cm
\noindent
{\bf A.2.4 Rankin-Selberg \& Zagier Method \newline
\hskip 1.13cm II: Semi-stable Case}

\smallskip
\noindent
The Fact G-I-J proves to be very important:
with this intrinsic relation between geo-arithmetical truncation and analytic 
truncation, instead of using geo-arithemtical method to study high rank zeta functions,  which is rather new and less developed,
we can equally use analytic technics and methods from trace formula, which is more systematic and rich, to help us. 
As an example, we here indicate how to evaluate the Eisenstein period $
\int_{\mathcal M_{\mathbb Q,r}[1]}E(\lambda;{\bf 1};g)\, dg.$ 

First, by Fact G-I-J, it is equal to
$\int_{G(\mathbb Z)\backslash G(\mathbb R)/SO(n)}\Lambda^0E(\lambda;{\bf 1};g)\, dg$. On the other hand, by Fact E$'$, we already know that when $T$ is sufficiently positive, $\int_{G(\mathbb Z)\backslash G(\mathbb R)/SO(n)}\Lambda^TE(\lambda;{\bf 1};g)\, dg$ can be evaluated. As such, then the only point here of course is to check whether the argument used for sufficiently positive $T$ are still valid when $T$ is taken to be $0$.

By examining the proof, to take care of the change from sufficiently positive $T$ to smaller $T$, say $T=0$, additional two main points must be checked. They are 

\noindent
(1) Fact G for smaller $T$. This now is replaced by Fact G-I-J. Cleared.

\noindent 
(2) The convergences of all integrations involved in the proof. This  is indeed a very serious one. In a  sense, modulo combinatorial technics,
establishing various convergences is really the technical heart of Arthur's trace formula
(in its preliminary form as stated in [Ar1-3]). Fortunately, we can justify these convergences when $T$ is smaller, in particular when $T=0$. Practically, this is carried out in two steps. First, for sufficiently positive $T$, we follow simply the original arguments in [Ar1-3] and [JLR]. Then for general $T\geq 0$, we use the fact that the difference for integral domains involved between sufficiently positive $T$ and rather small $T$, say, $T=0$, is only up to a certain suitable compact subset in a fundamental domain -- after all, over compact subsets,  integrability becomes trivial for smooth functions. In this way, we then arrive at the following
\vskip 0.30cm
\noindent
{\bf Fact E}$''$. ([W1,3]) {\it The Eisensetin period $\int_{\mathcal M_{\mathbb Q,r}[1]} E({\bf 1};\lambda;g)\, dg$ is given by
$$\int_{\mathcal M_{\mathbb Q,r}[1]} E({\bf 1};\lambda;g)\, dg=v\sum_{w\in W} \frac {1}{\prod_{\alpha\in \Delta_0} \langle w\lambda-\rho, \alpha^{\vee}\rangle}\cdot  M(w,\lambda)$$} 
\eject
\noindent
{\bf A.2.5 Intertwining Operator: Gindikin-Karpelevich Formula}

\smallskip
\noindent
To go further, we need write down also the intertwining operator $M(w,\lambda)$. This is now well known -- by the Gindikin-Karpelevich formula,
we have 

\noindent
{\bf Fact K.} (See e.g., [La3]) {\it For every split, semi-simple group $G$,
its associated intertwinging operator acting on constant function $\bold 1$ over the Borel is given by $$M(w,\lambda)=\prod_{\alpha>0, w\alpha<0} \frac {\xi\big(\langle\lambda, \alpha^{\vee}\rangle\big)}{\xi\big(\langle \lambda, \alpha^{\vee}\rangle+1\big)}.
$$
 Here $\xi(s)$ is the completed Riemann zeta  with $\Gamma$-factor, namely, $\xi(s)=\pi^{-\frac s2}\Gamma(\frac s2)\zeta(s)$ with $\zeta(s)=\sum_{n=1}^\infty n^{-s}$ the standard Riemann zeta function.}
\vskip 0.30cm
\noindent
{\bf A.2.6 Periods for $G$ over $\mathbb Q$: Weyl Symmetry}

\smallskip
\noindent
By Facts E$''$, K, for sufficiently positive $T$, the associated Eisenstein period $$\omega_{\mathbb Q}^{SL(n),T}(\lambda):=\int_{SL(n,\mathbb Z)\backslash SL(n,\mathbb R)/SO(n)} \wedge^T E({\bf 1};z_1,z_2,\dots, z_n;M)\, d\mu(M)$$ is given by the following
\vskip 0.30cm
\noindent
{\bf Fact L.} ([W1,3]) {\it  Up to a constant factor,
$$\omega_{\mathbb Q}^{SL(n),T}(\lambda)
=\sum_{w\in W} \frac {e^{\langle w\lambda-\rho, T\rangle}}{\prod_{\alpha\in \Delta_0} \langle w\lambda-\rho, \alpha^{\vee}\rangle}  \cdot \prod_{\alpha>0, w\alpha<0} \frac {\xi\big(\langle\lambda, \alpha^{\vee}\rangle\big)}{\xi\big(\langle \lambda, \alpha^{\vee}\rangle+1\big)}.
$$}

With this, by a close look at the right hand side, we conclude that
now we may take even $T=0$, even the right hand only makes sense for sufficiently positive $T\gg 0$.
This then leads to
\vskip 0.30cm
\noindent
{\bf Definition 1.} {\it The period for $G$ over $\mathbb Q$ is defined by
$$\omega_{\mathbb Q}^G(\lambda):=\sum_{w\in W}\Bigg(\frac{1}{\prod_{\alpha\in\Delta_0}\langle w\lambda-\rho,\alpha^\vee\rangle}\cdot\prod_{\alpha>0,w\alpha<0}\frac{\xi_{\mathbb Q}(\langle\lambda,\alpha^\vee\rangle)}
{\xi_{\mathbb Q}(\langle\lambda,\alpha^\vee\rangle+1)}\Bigg),\ \ \mathrm{Re}\,\lambda\in\mathcal C^+$$ where $\mathcal C^+$ denotes the standard positive Weyl chamber of $\frak a_0$, the space of 
characters associated to $(G,B)$, and 
$\xi_{\mathbb Q}(s)$ the completed Riemann zeta function.}

Certainly this is exact the definition 1 in the main text. As such, the most notable point in this definition is the huge symmetry created by the Weyl group.

\subsection{New Zetas for $SL(n)/\mathbb Q$}

\noindent
{\bf A.3.1 Epstein, Koecher, Siegel Zetas and Siegel-Eisenstein Series}

\smallskip
\noindent
The reason why we care about Eisenstein periods $\int_{\mathcal M_{\mathbb Q,r}[1]} E(\lambda;{\bf 1};g)\, dg$, which are of several variables, is that this period can be evaluated and that Epstein zetas $E(\Lambda^g,s)$ appeared in the study of high rank zetas are residues of Eisenstein series
$E({\bf 1};\lambda;g)$:  $$\xi_{\bf Q,r}(s)=\frac{r}{2}\int_{{\mathcal M}_{{\mathbb Q},r}[1]} \hat E(\Lambda^g,\frac{r}{2}s)\, dg$$ where $\hat E(\Lambda^g,s)=\pi^{-s}\Gamma(s)\cdot E(\Lambda^g,s)$. To explain this, we go as follows.

Let $\frak R:=\{\mathrm{diag}(\pm 1,\dots,\pm 1)\}\backslash 
SL(n,\mathbb Z)$ and $\frak Q_r$ the standard parabolic subgroup 
associated to the partition $n=r+1+1+\cdots+1$, that is, the parabolic subgroup $P_{r,1,\dots,1}$ consisting of matrices in $SL(n,\mathbb Z)$
of the form $
\begin{pmatrix}H&&&\\ &1&*&\\ &0&\ddots&\\ &&&1\end{pmatrix}$ with $H=H^{(r)}, |H|=1$.
Define the associated {\it Siegel zeta functions} by $$\xi_r^*(Y;s_r,\dots,s_{n-1}):=\sum_{N\in \frak Q_r\backslash \frak R}\prod_{v=r}^{n-1}|Y[N]_v|^{-s_v}$$ for all $1\leq r\leq n-1$, where, as usual, $Y[N]:=N^t\cdot Y\cdot N$ and for a size $n$ matrix $A=(a_{ij})_{i,j=1}^n$, $A_v$ denotes the matrix $A_v=(a_{ij})_{i,j=1}^v$ .
Then, from [D], we have the following

\noindent
{\bf Lemma 1.} ([D]) {\it There exists a constant $c$ depending only on $r$ such that} $$\mathrm{Res}_{s_r=\frac{r+1}{2}}\xi_r^*(Y;s_r,\dots,s_{n-1})=c_r\cdot \xi_{r+1}^*(Y;s_{r+1}+\frac{r}{2},s_{r+2},\dots,s_{n-1}).$$ (Please correct a misprint in [D] for this formula.)
Consequently, taking $r=1$ and repeating this process, we obtain the following
$$\mathrm{Res}_{s_{n-1}=1}\cdots \mathrm{Res}_{s_{2}=1}\mathrm{Res}_{s_{1}=1}\Big(\xi_1^*(Y;s_1,s_2,\dots, s_{n-1})\Big)=|Y|^{-\frac{n-1}{2}},$$ up to a constant factor.

Similarly, for the standard parabolic group $P=P_{n_1,n_2,\dots,n_q}$ corresponding to the partition $n=n_1+n_2+\cdots+n_q$,
define the associated {\it Siegel's Eisenstein series}  by
$$\begin{aligned}E_{P}({\bf s}|Y):=&E_{n_1,n_2,\dots,n_q}({\bf s}|Y)\\
:=&\prod_{(A_j*)=A\in\Gamma_n/P, A_j\in\mathbb Z^{n\times N_j}}\prod_{j=1}^q|Y[A_j]|^{-s_j},\qquad \mathrm{Re}s_j>\frac{n_j+n_{j+1}}{2}\end{aligned}$$ where ${\bf s}=(s_1,s_2,\dots,s_q)$ and $N_j=n_1+n_2+\cdots+n_j$.
Define also {\it Koecher's zeta function} by $$Z_{m,n-m}(X,s):=\sum_{A\in\mathbb Z^{n\times m}/GL(m,\mathbb Z),\mathrm{rk}A=m}|X[A]|^{-s},\qquad\mathrm{Re}(s)>\frac{n}{2}.$$ 

\noindent
{\bf Lemma 2.} (See e.g. [Te]) (1) {\it $E_{n_1,n_2,\dots,n_q}({\bf s}|Y)$ and $Z_{m,n-m}(X,s)$ are well-defined in the above indicated regions and admit meromorphic continuation to the whole parameter spaces;} and 

\noindent
(2) {\it They satisfy the following relations (see e.g. [Te]):
$$|Y|^{-s}\cdot E_{n-1,1}({\bf 1};s|Y^{-1})= E_{1,n-1}({\bf 1};s|Y)
=Z_{1,n-1}(Y;s)/Z_{1,0}(I;s)$$ and $$Z_{n,0}(X,s)=|X|^{-s}\cdot\prod_{j=0}^{n-1}\zeta(2s-j).$$}

In parallel, for a positive definite $Y$ with $|Y|=1$ and ${\bf s}=(s_1,s_2,\dots,s_n)$, introduce as usual the power function $$p_{-{\bf s}}(Y):=\prod_{j=1}^n|Y_j|^{-s_j}.$$ Then
the associated Siegel's Eisenstein series for the Borel $B=P_{1,1,\dots,1}$ is defined as  
$$E_{(n)}({\bf s}|Y):=\sum_{\gamma\in \Gamma_n/P_{1,1,\dots,1}}p_{-{\bf s}}(Y),\qquad\mathrm{Re}s_j>1, j=1,2,\dots, n-1.$$

\noindent
{\bf Lemma 3.} {\it We have
$$\xi_1^*(Y;s_1,s_2,\dots,s_{n-1})=E_{(n)}(s_1,s_2,\dots,s_n|Y),$$ and
$$E_{(n)}({\bf s}|Y^{-1})=E_{(n)}({\bf s}^*|Y)$$ where ${\bf s}^*:=(s_{n-1},s_{n-2},\dots,s_2,s_1,-(s_1+s_2+\dots+s_n))$.
Consequently, $$\xi_1^*(Y^{-1};t_1,t_2,\dots, t_{n-1})=
\xi_1^*(Y;t_{n-1},\dots,t_2, t_{1}).$$}

Thus, in particular, for the Siegel Eisenstein series corresponding to the
maximal parabolic subgroup $P_{n-1,1}$, i.e., for $$\begin{aligned}E_{n-1,1}(s_1,s_2|Y):=&E_{n-1,1}({\bf 1};s_1,s_2|Y)\\
:=&\sum_{(A_1*)=A\in \Gamma_n/P_{n-1,1}, A_1\in\mathbb Z^{n\times(n-1)}}|Y[A_1]|^{-s_1}|Y[A]|^{-s_2},\end{aligned}$$ 
 we have, by Lemma 3,
$$\begin{aligned}E_{n-1,1}(s,t|Y)=&|Y|^{-t}\cdot\sum_{(A_1*)=A\in\Gamma_n/P_{n-1,1}}|Y[A_1]|^{-s}\\
=&|Y|^{-t}\cdot\sum_{A\in\Gamma_n/P_{n-1,1}}|Y[A]_{n-1}|^{-s}=\xi^*_{n-1}(Y,s).\end{aligned}$$ Here, we used the fact that the group involved is $SL(n)$.

\bigskip
Consequently, by Lemmas 1 and 2, we obtain the following

\noindent
{\bf Fact M.} (1)  {\it $\xi^*_{n-1}(Y;s)$ and $E(\Lambda;s)$ are related by $$\xi_{n-1}^*(Y^{-1};s)=\frac{1}{\zeta(2s)}\cdot\sum_{{\bf x}\in\mathbb Z^n\backslash\{0\}}|Y[{\bf x}]|^{-s}=\frac{1}{\zeta(2s)}\cdot E\Big(\Lambda(g);\frac{s}{n/2}\Big)$$
where $Y:=g^t\cdot g$ and $\Lambda$ denotes the lattice $(\mathbb Z^n;\rho(g))$ with the metric $\rho(g)$ on $\mathbb R^n$ induced by the positive definite matrix $Y=g^t\cdot g$;} and

\noindent
(2) $\displaystyle{\xi_1^*(Y;s_1,s_2,\dots,s_{n-1})=E_{(n)}(s_1,s_2,\dots,s_n|Y).}$

\smallskip
\noindent
{\it In particular,}

\smallskip
\noindent
 $E(\Lambda(g);s)=\mathrm{Res}_{t_{n-2}=1,\,t_{n-3}=1,\,\dots,\, t_{2}=1,\,t_{1}=1}\xi_1^*(Y;ns-\frac{n-2}{2},t_{n-2},t_{n-3},\dots, t_2,t_1).$

\vskip 0.30cm
\noindent
{\bf A.3.2 Langlands' Eisenstein and Siegel's Eisenstein}

\smallskip
\noindent
To apply Fact E$''$ directly, we still need write Siegel's Eisenstein series introduced using classical language in terms of Langlands' Eisenstein series introduced using a language which is more convenient for theoretical purpose.
The point of course is about the power function $p$ and the function $m_B$. For this, write a positive definite $Y$ (with $|Y|=1$) as $Y={\bf a}[{\bf n}]$ with ${\bf a}=\mathrm{diag}(a_1,a_2,\dots,a_n)$ and ${\bf n}$ upper triangular unipotent (with diagonal entries 1).
Then $a_i=|Y_i|/|Y_{i-1}|, i=1,2,\dots,n$.
Consequently, by definition,
$$\begin{aligned}&p_{-{\bf s}}(Y)=\prod_{j=1}^n|Y_j|^{-s_j}\\
=&|Y_1|^{-(s_1+s_2+\cdots+s_n)}\Big(|Y_2|/|Y_1|\Big)^{-(s_2+s_3+\cdots+s_n)}\\
&\qquad\cdots
\Big(|Y_{n-1}|/|Y_{n-2}|\Big)^{-(s_{n-1}+s_n)}\cdot \Big(|Y_{n}|/|Y_{n-1}|\Big)^{-s_n}\\
=&a_1^{-(s_1+s_2+\cdots+s_n)}a_2^{-(s_2+s_3+\cdots+s_n)}\cdots
a_{n-1}^{-(s_{n-1}+s_n)}\cdot \Big(a_1a_2\cdots a_{n-1}\Big)^{-s_n}\\
=&a_1^{-(s_1+s_2+\cdots+s_{n-1})}a_2^{-(s_2+s_3+\cdots+s_{n-1})}\cdots
a_{n-1}^{-s_{n-1}}\end{aligned}$$ since $\prod_{j=1}^na_j=|Y|=1$.

On the other hand, if $Y=g^tg$ with $T(g)=\mathrm{diag}(t_1,t_2,\dots,t_n)$, then we have $a_j=t_j^2$ and
$$m_B(g)^{\lambda+\rho_B}=T(g)^{\lambda+\rho_B}$$
where as usual, we let $\lambda=(z_1,z_2,\dots,z_n)\in\mathbb C^n,\ \sum_{j=1}^nz_j=0$ so that $$\rho=\rho_B=\Big(\frac{n-1}{2},\frac{n-1}{2}-1,\dots, 1-\frac{n-1}{2},-\frac{n-1}{2}\Big).$$
Hence, by a direct calculation, we get $$\begin{aligned}&m_B(g)^{\lambda+\rho_B}\\
=&t_1^{-[(n-1)+(2z_1+z_2+\cdots+z_{n-1})]}\cdot t_2^{-[(n-2)+(z_1+2z_2+\cdots+z_{n-1})]}\cdots t_{n-1}^{-[1+(z_1+z_2+\cdots+2z_{n-1})]}.\end{aligned}$$ 

Recall also that the Langlands Eisenstein series associated to the constant function ${\bf 1}$ on the Borel $B=P_{1,1,\dots,1}$
(related to  $SL(n)/B$) is given by $$E({\bf 1};\lambda)(g):=\sum_{\gamma\in SL(n,\mathbb Z)/P_{1,1,\dots,1}=B}m_B(\delta g)^{\lambda+\rho_B}.$$ So if we make the variable transformation
from $\lambda$ to ${\bf s}$ by
$$\begin{cases}2s_1=&1+(z_1-z_2)\\
2s_2=&1+(z_2-z_3)\\
\cdots&\cdots\\
2s_{n-1}=&1+(z_{n-1}-z_{n})\end{cases}$$
Then we arrive at the

\noindent
{\bf Fact M}$'$. (1) {\it $E({\bf 1};\lambda)(g)=E_{(n)}({\bf s}|Y)$, \newline where 
$\lambda=(z_1,z_2,\dots,z_n)$ with $\sum_{j=1}^nz_j=0$ and ${\bf s}=(s_1,s_2,\dots,s_{n-1})$ satisfying
$$\begin{cases}2s_1=&1+(z_1-z_2)\\
2s_2=&1+(z_2-z_3)\\
\cdots&\cdots\\
2s_{n-1}=&1+(z_{n-1}-z_{n}).\end{cases}$$}

\noindent
(2) {\it Introduce the variable $s$ via $2ns-n+1=z_1-z_2,$, we have
 $$E\Big(\Lambda(g);s\Big)=\mathrm{Res}_{z_2-z_3=1}\mathrm{Res}_{z_3-z_4=1}\cdots
\mathrm{Res}_{z_{n-1}-z_n=1}E({\bf 1};z_1,z_2,\dots,z_n)(g).$$}
\vskip 0.20cm
\noindent
{\bf A.3.3 New Zetas: Genuine but Different}

\smallskip
\noindent
Recall that, by Fact D, high rank zetas are given by
$$\xi_{\mathbb Q,r}(s)=\int_{\mathcal M_{\mathbb Q,r}[1]}
\widehat E(\Lambda;\frac{r}{2}s)\,d\mu_0(\Lambda).$$ Thus by Facts G-I-J and M, to offer a close formula, it suffices to  evaluate the integration
$$\int_{\frak F(0)}\mathrm{Res}_{z_2-z_3=1}\cdots \mathrm{Res}_{z_{3}-z_4=1}\cdots \mathrm{Res}_{z_{r-1}-z_r=1}\Big(E({\bf 1};z_1,z_2,\dots, z_r)(g)\Big)\,d\mu(g).$$ Thus, if we were able to freely make an interchange between

\noindent
(i) the operation of taking integration $\int_{\frak F(0)}$ and 

\noindent
(ii) the operation of taking residues $\mathrm{Res}_{z_2-z_3=1}\cdots \mathrm{Res}_{z_{3}-z_4=1}\cdots \mathrm{Res}_{z_{r-1}-z_r=1}$,

\noindent
it would be sufficient for us to evaluate 
$$\mathrm{Res}_{z_2-z_3=1}\cdots \mathrm{Res}_{z_{3}-z_4=1}\cdots \mathrm{Res}_{z_{r-1}-z_r=1}\Big(\int_{\frak F(0)}E({\bf 1};z_1,z_2,\dots, z_r)(g)\,d\mu(g)\Big),$$ or better, to evaluate the expression
$$\begin{aligned}\mathrm{Res}_{z_2-z_3=1}&\cdots \mathrm{Res}_{z_{3}-z_4=1}\cdots \mathrm{Res}_{z_{r-1}-z_r=1}\\
&\Bigg(\sum_{w\in \Omega} \frac {1}{\prod_{\alpha\in \Delta_0} \langle w\lambda-\rho, \alpha^{\vee}\rangle}  \cdot \prod_{\alpha>0, w\alpha<0} \frac {\xi\big(\langle\lambda, \alpha^{\vee}\rangle\big)}{\xi\big(\langle \lambda, \alpha^{\vee}\rangle+1\big)}\Bigg)\end{aligned}$$ where $\lambda=(z_1,z_2,\dots,z_r)$ with $z_1+z_2+\cdots+z_r=0$, since by Fact G,
$$\begin{aligned}\int_{\frak F(T)}&E({\bf 1};z_1,z_2,\dots, z_r)(g)\,d\mu(g)\\
=&\int_{SL(r,\mathbb Z)\backslash SL(r,\mathbb R)/SO(r)}\Lambda^TE({\bf 1};z_1,z_2,\dots, z_r)(g)\,d\mu(g)\\
=&\sum_{w\in \Omega} \frac {e^{\langle w\lambda-\rho, T\rangle}}{\prod_{\alpha\in \Delta_0} \langle w\lambda-\rho, \alpha^{\vee}\rangle}  \cdot \prod_{\alpha>0, w\alpha<0} \frac {\xi\big(\langle\lambda, \alpha^{\vee}\rangle\big)}{\xi\big(\langle \lambda, \alpha^{\vee}\rangle+1\big)}\end{aligned}$$
by Fact L.

Unfortunately, this interchange of orders of two operations is not allowed in general. As examples, 
one can observe this by working on $SL(n)$ and by comparing the poles for the resulting expressions. (For details, see the remark at the end of A.3.4 below.)
\vskip 0.30cm
On the other hand, even with the existence of such discrepancies, the function 
$$\begin{aligned}\mathrm{Res}_{z_2-z_3=1}&\cdots \mathrm{Res}_{z_{3}-z_4=1}\cdots \mathrm{Res}_{z_{r-1}-z_r=1}\\
&\Bigg(\sum_{w\in W} \frac {1}{\prod_{\alpha\in \Delta_0} \langle w\lambda-\rho, \alpha^{\vee}\rangle}  \cdot \prod_{\alpha>0, w\alpha<0} \frac {\xi\big(\langle\lambda, \alpha^{\vee}\rangle\big)}{\xi\big(\langle \lambda, \alpha^{\vee}\rangle+1\big)}\Bigg)\end{aligned}$$ proves to be extremely {\it natural and nice}. This then leads to 
\vskip 0.30cm
\noindent
{\bf Definition 2.} The single variable period $Z_{\mathbb Q}^{SL(r)}(z_1)$ associated to $SL(r)$ over $\mathbb Q$ is defined by
$$\begin{aligned}Z_{\mathbb Q}^{SL(r)}(z_1):=&\mathrm{Res}_{z_2-z_3=1}\cdots \mathrm{Res}_{z_{3}-z_4=1}\cdots \mathrm{Res}_{z_{r-1}-z_r=1}\\
&\Bigg(\sum_{w\in W} \frac {1}{\prod_{\alpha\in \Delta_0} \langle w\lambda-\rho, \alpha^{\vee}\rangle}  \cdot \prod_{\alpha>0, w\alpha<0} \frac {\xi\big(\langle\lambda, \alpha^{\vee}\rangle\big)}{\xi\big(\langle \lambda, \alpha^{\vee}\rangle+1\big)}\Bigg),\end{aligned}$$ where $\lambda=(z_1,z_2,\dots,z_r)$ with $z_1+z_2+\cdots+z_r=0$.

Clearly, there are some factors $\xi(ax+b)$'s  
left in the denominator even after all cancelations made. To clear them, we make the following observations:

\noindent
(i) there is a minimal integer $I(SL(r))$ and finitely many factors $$\xi\Big(a_1^{SL(r)}z_1+b_1^{SL(r)}\Big),\ \xi\Big(a_2^{SL(r)}z_1+b_2^{SL(r)}\Big),\ \cdots,\ 
\xi\Big(a_{I(SL(r))}^{SL(r)}\lambda_P+b_{I(SL(r))}^{SL(r)}\Big),$$
such that {\it the product $\Big[\prod_{i=1}^{I(SL(r))}\xi\Big(a_i^{SL(r)}z_1+b_i^{SL(r)}\Big)\Big]\cdot Z_{\mathbb Q}^{SL(r)}(z_1)$ admits only finitely many singularities.}

\noindent
(ii) there is a minimal integer $J(SL(r))$ and finitely many factors 
$$\xi\Big(c_1^{SL(r)}\Big),\ \xi\Big(c_2^{SL(r)}\Big),\
\cdots,\
\xi\Big(c_{J(SL(r))}^{SL(r)}\Big),$$ 
such that
{\it there are no factors of special $\xi$ values appearing at the denominators in  the product  $\Big[\prod_{i=1}^{J(SL(r))}\xi\Big(c_i^{SL(r)}\Big)\Big]\cdot Z_{\mathbb Q}^{SL(r)}(z_1).$} 
\vskip 0.30cm
\noindent
{\bf Definition 3.} The {\it zeta function $\xi_{\mathbb Q;o}^{SL(r)}$ for $SL(r)$ over $\mathbb Q$} is defined by
$$\begin{aligned}\xi_{\mathbb Q;o}^{SL(r)}\Big(s\Big)
:=&\Bigg[\prod_{i=1}^{I(SL(r))}\xi\Big(a_i^{SL(r)}s+b_i^{SL(r)}\Big)\cdot\prod_{j=1}^{J(SL(r))}\xi\Big(c_j^{SL(r)}\Big)\Bigg]\cdot Z_{\mathbb Q}^{SL(r)}\Big(s\Big),\\
&\hskip 7.5cm\mathrm{Re}\,s\gg 0\end{aligned}
$$

Clearly, Definitions 2 and 3 here are special cases of Definitions 2 and 3 in the main text. In fact here implicitly the maximal parabolic subgroup $P_{r-1,1}$ is used.

\vskip 0.30cm
\noindent
{\it Remark.} We remind the reader that the version with parameter $T$ is in fact also very important. In rank two case, one can show that for $T$ non-negative, the associated period also satisfies the functional equation and the  RH. For general cases, the structure is more complicated on one hand,
and beautiful on the other:
Say the functional equation for $\xi_{\mathbb Q}^{SL(n)/P_{m,n-m};T}$ is related with a different function $\xi_{\mathbb Q}^{SL(n)/P_{n-m,m};T}$ (for a different maximal parabolic subgroup), based on
another type of symmetry between $E_{m,n-m}$ for $Y$ and $E_{n-m,n}$ for $Y^{-1}$ stated above (for classical Siegel Eisenstein series).
However, when $T=0$, $\xi_{\mathbb Q}^{SL(n)/P_{m,n-m};0}$ is essentially the function $\xi_{\mathbb Q}^{SL(n)/P_{n-m,m};0}$.
All this then leads to the functional equation for $\xi_{\mathbb Q}^{SL(n)/P_{m,n-m}}(s)$.
\vskip 0.30cm
\noindent
{\bf A.3.4 Functional Equation and the Riemann Hypothesis}

Just as high rank zetas, we certainly expect that these new zetas
introduced in the previous subsection
satisfy the functional equation and an analogue of the Riemann Hypothesis.
For this we have the following
\vskip 0.30cm
\noindent
{\bf Conjecture.} ({\bf Functional Equation}) {\it There exists a constant $c_{SL(r)}$ depending on $r$ only such that
$$\xi_{\mathbb Q;o}^{SL(r)}\Big(c_{SL(r)}-s\Big)=\xi_{\mathbb Q;o}^{SL(r)}\Big(s\Big).$$}

\vskip 0.30cm
To make the functional equation canonical, i.e., reflecting the standard symmetry $s\leftrightarrow 1-s$ for the standard functional equation, we make the following  normalization.

\noindent
{\bf Definition 3}$'$  The {\it zeta function $\xi_{\mathbb Q}^{SL(r)}\Big(s\Big)$ for $SL(r)$ over $\mathbb Q$} is defned by
$$\xi_{SL(r);\mathbb Q}\Big(s\Big):=\xi_{\mathbb Q;o}^{SL(r)}\Big(s+\frac{c_{SL(r)}-1}{2}\Big)$$

As such then we have the following
\vskip 0.30cm
\noindent
{\bf Conjecture$'$.} ({\bf Functional Equation}) {\it $\xi_{SL(r);\mathbb Q}(1-s)=\xi_{SL(r);\mathbb Q}(s).$}

\vskip 0.30cm
The most remarkable property shared by all these newly introduced zetas
is the following Zeta Fact about the uniformity of their zeros. That is to say, we expect the following
\vskip 0.20cm
\noindent
{\bf The Riemann Hypothesis\,$^{G/P}_{\mathbb Q}$.}

\noindent
$$All\ zeros\ of\ the\ zeta\ function\ \xi_{SL(r);\mathbb Q}(s)\ lie\ on\ the\ central\ line\
\mathrm{Re}\,s=\displaystyle{\frac{1}{2}}$$

\noindent
{\it Remark.} Examples with $SL(3)$ shows that $\xi_{\mathbb Q,r}(s)$ is not the same as $\xi_{SL(r);\mathbb Q}(s)$ in general.
(In fact, while $\xi_{{\mathbb Q},3}(s)$ has only two singularities at $s=0,1$,
 $\xi_{{\rm SL}(3);{\mathbb Q}}(s)$ has four singularities at $s=0,\,\frac{1}{3},\,\frac{2}{3},\,1$, by the precise formula listed in Appendix B.)
 We detect such a discrepancy only quite later, indeed, not even until the first version of this paper was written on Dec. 4, 2007: so quite some time, we wrongly believed that $\xi_{\mathbb Q,r}(s)$ are $\xi_{SL(r);\mathbb Q}(s)$ are the same.
\vskip 0.30cm
After making these conjectures, or more correctly, after making the RH for high rank zetas open, (with the mistake mentioned in the remark above,) we felt that more examples should be provided at least numerically. This then led to the problem of finding precise expressions for \lq$\xi_{\mathbb Q,r}(s)$' with $r=4,5$. For quite some time, we could not make an advance. However, this was changed with the short visit of Henry Kim in the summer of 2007, who brought us the paper of Diehl [D]. With [D], being compatible with our old approach for rank 3 zeta by taking residues in [W3],
we were afterwards able to see the structure for the zetas, modulo a few mistakes.
Accordingly, we did some painful calculations:

\noindent
a) For rank 4, totally $24\times 6=144$ cases were discussed, from which we obtained the final zeta consisting of 12 terms;

\noindent
b) For rank 5, totally $120\times 10=1200$ cases were discussed, from which we obtained the final zeta consisting of 28 terms.

For details, see Appendix B: Examples.
Based on all these calculations, we are able to exposes the following
\vskip 0.30cm
\noindent
{\bf Fact N.} (1) ({\bf Functional Equation}$_{\leq 5}$) $$\xi_{SL(r);\mathbb Q}(1-s)=\xi_{SL(r);\mathbb Q}(s)\qquad {when}\ r=2,\,3,\,4,\,5;$$ 

\noindent
(2) ([LS], [S]) ({\bf Riemann Hypothesis}$_{SL(2,\,3);\mathbb Q}$)

{\it All zeros of $\xi_{SL(2);\mathbb Q}(s)$ and $\xi_{SL(3);\mathbb Q}(s)$ lie on the 
central line $\mathrm{Re}(s)=\frac{1}{2}.$}

\subsection{Zetas for $(G,P)/\mathbb Q$}
\vskip 0.20cm
\noindent
{\bf A.4.1 From SL to Sp: Analytic Method Adopted \& Periods Chosen}

\smallskip
\noindent
For quite sometime, we want to use geo-arithmetic method to find an analogue of high rank zetas for other reductive groups. The first natural target is $Sp$. However, this proves to be very difficult, since for the completed theory, we should start with what might be called principal lattices associated to $Sp$ and establish all the $Sp$ properties corresponding to Facts listed above for $SL$.

Fortunately, for the purpose of finding corresponding zeta functions  $\xi_{Sp(2n);\mathbb Q}(s)$ for $Sp$,
with our success for $SL$ discussed above and the paper of Diehl [D], which in fact deals with $Sp$ instead of $SL$, we realize that instead of the approach using geo-arithmetic method, an alternative way using pure analytic methods is sufficient.  This goes as follows.

Let $G=Sp(2n)$ with $G(\mathbb R)=Sp(2n,\mathbb R)$ the symplectic group of degree $n$ over $\mathbb R$.
For any $Z\in \frak S=\frak S_n$, the Siegel upper half space of rank $n$,  write $Z=X+\sqrt{-1}Y$ according to its real and imaginary parts. By definition, $Y=\mathrm{Im}\,Z>0$ and $Z^t=Z$ is symmetric. For an $M=\begin{pmatrix} A&B\\ C&D\end{pmatrix}\in Sp(2n,\mathbb R)$, as usual, set
$M\langle Z\rangle:=(AZ+B)\cdot(CZ+D)^{-1}$ and write $Y(M):=\mathrm{Im} M\langle Z\rangle$. Note that the action is transitive and the stablizer in $Sp(2n,\mathbb R)$ for $\sqrt {-1}I$ is given by $Sp(2n,\mathbb R)\cap SO(2n)$.
Consequently, we obtain a natural isomorphism $Sp(2n,\mathbb R)/SO(2n)\cap Sp(2n,\mathbb R)\simeq {\frak S}_n$.

Introduce  also  $\Gamma_n:=\big\{\mathrm{diag}(\pm 1,\pm1,\cdots,\pm 1)\big\}\backslash Sp(2n,\mathbb Z)$ the Siegel modular group,
and $\frak B=\frak P_n:=\Big\{\begin{pmatrix} *&*\\ 0&*\end{pmatrix}\in \Gamma\Big\}$ the associated maximal parabolic subgroup.

Fix $Z\in \frak S$, define then the associated {\it Siegel-Maa$\beta$} Eisenstein series, or better, the {\it Siegel-Epstein zeta function} by
$$E_n(Z;s):=\sum_{\gamma\in\frak B\backslash\Gamma}\frac{|Y|^{-s}}{\|CZ+D\|^{-2s}}.$$

Motivated by our study on high rank zetas associated to $SL(n)$, for sufficiently positive $T$,
we define a {\it principal period for $Sp(n)$ over $\mathbb Q$}  by 
$$\zeta_{Sp(n),\mathbb Q}^T(s):=\int_{\Gamma\backslash \frak S_n}\wedge^T E_n(Z;s)\,d\mu(Z).$$
This is then a function on $s$ depending also on the parameter $T$.
It is then an open problem whether we can  evaluate this expression at $T=0$ since the corresponding Fact G-I-J for $Sp$ is still missing. Assume that the answer to this is affirmative, then $$\zeta_{Sp(n),\mathbb Q}(s):=\zeta_{Sp(n),\mathbb Q}^0(s):=\zeta_{Sp(n),\mathbb Q}^T(s)|_{T=0}$$ may be viewed as an $Sp$-analogue of the high rank zeta functions, call it
the {\it principal zeta function for $Sp(n)$ over $\mathbb Q$}.

As for the case of $SL(n)$, it is, for the time being, {\it very difficult}, in fact, {\it quite impossible,}  to offer a precise formula for the Eisenstein period $\zeta_{Sp(n),\mathbb Q}^T(s)$. 
However, motivated by our study for $SL(n)$, we want to introduce an analogue for the new type of zeta functions $\xi_{SL(r);\mathbb Q}(s)$.
For this (a bit changed yet very meaningful) purpose, we make the following preparations.
\vskip 0.30cm
\noindent
a) {\bf Siegel Eisenstein series.} As usual, corresponding to the partition $n=r+1+1+\cdots+1$,
introduce the standard parabolic subgroup $\frak P_r:=\Big\{\begin{pmatrix} A&*\\ 0&B\end{pmatrix}\in\Gamma\Big\}
$ where $A=\begin{pmatrix}H^t&&&\\ &1&0&\\ &*&\ddots&\\ &&&1\end{pmatrix}, B=\begin{pmatrix}H^{-1}&&&\\ &1&*&\\ &0&\ddots&\\ &&&1\end{pmatrix}$ with $H=H^{(r)}, |H|=1$. Accordingly, define the associated {\it Siegel Eisenstein series} by $$E_r(Z;s_r,\dots,s_n):=\sum_{\gamma\in\frak P_r\backslash \Gamma}\prod_{v=r}^n|Y(\gamma)_v)|^{-s_v}.$$
It is known that these Siegel Eiesnetsin series are naturally related to the Siegel zeta functions associated to the
standard parabolic subgroup $\frak Q_r$ of $SL(n)$, used already in our study for zetas associated to $SL(n)$. Recall that, if $\frak R:=\{\mathrm{diag}(\pm 1,\dots,\pm 1)\}\backslash SL(n,\mathbb Z)$ and $\frak Q_r$ is the standard parabolic subgroup
associated to the partition $n=r+1+1+\cdots+1$, then the associated {\it Siegel zeta functions} are defined by $$\xi_r^*(Y;s_r,\dots,s_{n-1}):=\sum_{N\in \frak Q_r\backslash \frak R}\prod_{v=r}^{n-1}|Y[N]_v|^{-s_v}$$ for all $1\leq r\leq n-1$.

\vskip 0.30cm
\noindent
{\bf Lemma 1.} ([D]) We have

\noindent
(i) $$E_r(Z;s_r,\dots, s_n)=\sum_{\gamma\in \frak B\backslash\Gamma}|Y(\gamma)|^{-s_n}\cdot
\xi_r^*\Big(Y(\gamma);s_r,\dots,s_{n-1}\Big);$$

\noindent
(ii) {\it There exists a constant $c$ depending only on $r$ such that} $$\mathrm{Res}_{s_r=\frac{r+1}{2}}\xi_r^*(Y;s_r,\dots,s_{n-1})=c_r\xi_{r+1}^*(Y;s_{r+1}+\frac{r}{2},s_{r+2},\dots,s_{n-1}).$$ 
Consequently, 
$$\mathrm{Res}_{s_{n-1}=1}\cdots \mathrm{Res}_{s_{2}=1}\mathrm{Res}_{s_{1}=1}\Big(\xi_1^*(Y;s_1,s_2,\dots, s_{n-1})\Big)=|Y|^{-\frac{n-1}{2}}$$ up to a constant factor.
Therefore, up to constant factors,
$$\begin{aligned}&\mathrm{Res}_{s_{n-1}=1}\cdots \mathrm{Res}_{s_{2}=1}\mathrm{Res}_{s_{1}=1}E_r(Z;s_r,\dots, s_n)\\
=&\sum_{\gamma\in \frak B\backslash\Gamma}|Y(\gamma)|^{-s_n}\cdot
\mathrm{Res}_{s_{n-1}=1}\cdots \mathrm{Res}_{s_{2}=1}\mathrm{Res}_{s_{1}=1}\xi_r^*(Y(\gamma);s_r,\dots,s_{n-1})\\
=&\sum_{\gamma\in \frak B\backslash\Gamma}|Y(\gamma)|^{-s_n}\cdot
|Y(\gamma)|^{-\frac{n-1}{2}}=E_n(Z;s_n+\frac{n-1}{2}).\end{aligned}$$

\noindent
b) {\bf Siegel Eisenstein series and Langlands Eisenstein series.}
As for the case of $SL(n)$, we next write the classical Siegel Eisenstein series in terms of Langlands' language. This is given by the following formula:
Let $\lambda=(z_1,z_2,\dots,z_n)\in\frak a_0$, then by defintion, $${\bf a}^{\lambda}(Z)=\prod_{v=1}^na_v^{-z_v}\qquad \mathrm{with}\qquad a_v=|Y_v|/|Y_{v-1}|.$$ Thus, the so-called power function
$${\bf p}_{-\bf s}(Y):=\prod_{\mu=1}^n|Y_\mu|^{-s_\mu}$$ is given by
$$\begin{aligned}\prod_{\mu=1}^n|Y_\mu|^{-s_\mu}=&{\bf p}_{-\bf s}(Y)={\bf a}^{\lambda}(Y)
=\prod_{v=1}^na_v^{-z_v}\\
=&|Y_1|^{-z_1+z_2}|Y_2|^{-z_2+2_3}\cdots|Y_{n-1}|^{-z_{n-1}+z_n}|Y_n|^{-z_n}.\end{aligned}$$
That is to say, we need  make the following change of variables $$s_1=z_1-z_2,
s_2=z_2-z_3,
\dots,
s_{n-1}=z_{n-1}-z_n,
s_n=z_n.$$
Consequently, we obtain the following 

\noindent
{\bf Fact M}$''$. (1) $E({\bf 1};\lambda;Y)=E_1(Z;s_1,s_2,\dots,s_n)$, and 

\noindent
(2) {\it Up to a suitable constant factor, 
$$\begin{aligned}&E_n(Z,z_n+\frac{n-1}{2})\\
=&\mathrm{Res}_{z_{n-1}-z_n=1}\cdots \mathrm{Res}_{z_{2}-z_3=1}\mathrm{Res}_{z_{1}-z_2=1}E({\bf 1};z_1,z_2,\dots, z_n;Y).\end{aligned}$$}

In particular, when $n=2$, i.e, for $Sp(4)$, we have 
$$\mathrm{Res}_{z_{1}-s=1}E({\bf 1};z_1,s;Y)=E_n(Z,s+\frac{1}{2}).$$

\noindent
c) {\bf The Sigel-Maa$\beta$-Eisenstein period.}
Note that the constant function one on the Borel is {\it cuspidal}, by the result of [JLR] cited above, and using the corresponding Gindikin-Karpelevich formula
for the associated intertwining operator,
 we have the following:
\vskip 0.30cm
\noindent
{\bf Fact E}$^{(3)}$. {\it Up to a constant factor,} $$\begin{aligned}&\int_{Sp(n,\mathbb Z)\backslash \frak S_n} \wedge^T E({\bf 1};\lambda;M)\, d\mu(M)\\
=&\sum_{w\in W} \frac {e^{\langle w\lambda-\rho, T\rangle}}{\prod_{\alpha\in \Delta_0} \langle w\lambda-\rho, \alpha^{\vee}\rangle}  \cdot \prod_{\alpha>0, w\alpha<0} \frac {\xi\big(\langle\lambda, \alpha^{\vee}\rangle\big)}{\xi\big(\langle \lambda, \alpha^{\vee}\rangle+1\big)}.\end{aligned}
$$ 

With all this, we are now ready to introduce our new zeta for $Sp(2n)$:
first define  ({\it not-yet-normalized}) {\it zeta} as the residue
$$\begin{aligned}&\mathrm{Res}_{z_{n-1}-z_n=1,\cdots,z_2-z_3=1, z_2-z_1=1}\\
&\qquad\sum_{w\in W} \frac {1}{\prod_{\alpha\in \Delta_0} \langle w\lambda-\rho, \alpha^{\vee}\rangle}  \cdot \prod_{\alpha>0, w\alpha<0} \frac {\xi\big(\langle\lambda, \alpha^{\vee}\rangle\big)}{\xi\big(\langle \lambda, \alpha^{\vee}\rangle+1\big)},\end{aligned}$$ since $\langle \rho,\alpha^\vee\rangle=1$ for all $\alpha\in\Delta_0$, where $\lambda=(z_1,z_2,\dots,z_n)\in\frak a_0$, corresponding to Definition 2;  then, make  certain normalizations following Definition 3. As such, we finally obtain a new series natural zetas $\xi_{Sp(2n),\mathbb Q}(s)$ for $Sp(2n)$ over $\mathbb Q$, which in fact coincide with $\xi_{\mathbb Q}^{Sp(2n)/\frak P_n}(s)$ defined in the main text.

As concrete examples, we worked out all the details for $n=2$. Similarly, we have the functional equation $$\xi_{Sp(4),\mathbb Q}(1-s)=\xi_{Sp(4),\mathbb Q}(s).$$   For details, see Appendix B below.

In summary, what we have done for $Sp$ is as follows:

\noindent
(i) First, motivated by our study for high rank zeta functions associated to $SL(n)$, we introduce a principal zeta for $Sp(2n)$ by evaluating the integration 
$$\int_{Sp(2n,\mathbb Z)\backslash\frak S_n}\wedge^TE_n(Z;s)\,d\mu(Z)$$ at $T=0$: in assuming that Fact G-I-J for $Sp$ can be established, even in the integration $T$ is supposed to be sufficiently positive, an evaluation at $T=0$ is allowed;

\noindent
(ii) By b), we know that, up to constant factors,
$$E_n(Z;z_n)=\mathrm{Res}_{z_{n-1}-z_n=1,\cdots,z_2-z_3=1, z_2-z_1=1}E({\bf 1};z_1,z_2,\dots, z_{n-1},z_{n}+\frac{n-1}{2};Y).$$ So it suffices to evaluate
$$\begin{aligned}\int_{Sp(2n,\mathbb Z)\backslash\frak S_n}&\mathrm{Res}_{z_{n-1}-z_n=1,\cdots,z_2-z_3=1, z_2-z_1=1}\\
&\qquad\Big(\wedge^TE({\bf 1};z_1,z_2,\dots, z_{n-1},z_{n}+\frac{n-1}{2};Y)\Big)\,d\mu(Y);\end{aligned}$$

\noindent
(iii) Even an interchange of $\int_{Sp(2n,\mathbb Z)\backslash\frak S_n}$ and $\mathrm{Res}_{z_{n-1}-z_n=1,\cdots,z_2-z_3=1, z_2-z_1=1}$ is not allowed, we, motivated by our success for $SL(n)$, still decide to study the period
$$\begin{aligned}&\mathrm{Res}_{z_{n-1}-z_n=1,\cdots,z_2-z_3=1, z_2-z_1=1}\\
&\qquad\int_{Sp(2n,\mathbb Z)\backslash\frak S_n}\Big(\wedge^TE({\bf 1};z_1,z_2,\dots, z_{n-1},z_{n}+\frac{n-1}{2};Y)\Big)\,d\mu(Y);\end{aligned}$$

\noindent
(iv) Now by c), for sufficiently positive $T$,
$\int_{Sp(2n,\mathbb Z)\backslash\frak S_n}\wedge^TE({\bf 1};\lambda;Y)\,d\mu(Y)$ is simply $$\sum_{w\in W} \frac {e^{\langle w\lambda-\rho, T\rangle}}{\prod_{\alpha\in \Delta_0} \langle w\lambda-\rho, \alpha^{\vee}\rangle}  \cdot \prod_{\alpha>0, w\alpha<0} \frac {\xi\big(\langle\lambda, \alpha^{\vee}\rangle\big)}{\xi\big(\langle \lambda, \alpha^{\vee}\rangle+1\big)}.$$ 

\noindent
(v) Evaluate the latest period at $T=0$ using the expression appeared in the right hand side and further take the residue.
This then leads to the not yet normalized new zeta function for $Sp(2n)$ over $\mathbb Q$: 
$$\begin{aligned}&\mathrm{Res}_{z_{n-1}-z_n=1,\cdots,z_2-z_3=1, z_2-z_1=1}\\
&\qquad\Big(\sum_{w\in W} \frac {1}{\prod_{\alpha\in \Delta_0} \langle w\lambda-\rho, \alpha^{\vee}\rangle}  \cdot \prod_{\alpha>0, w\alpha<0} \frac {\xi\big(\langle\lambda, \alpha^{\vee}\rangle\big)}{\xi\big(\langle \lambda, \alpha^{\vee}\rangle+1\big)}\Big).\end{aligned}$$

\noindent
(vi) Suitably normalized, we obtain a new type of zeta function $\xi_{Sp(2n)\mathbb Q}(s)$ for which we have the following
\vskip 0.30cm
\noindent
{\bf Conjecture.} (1) ({\bf Functional Equation}) {\it $\xi_{Sp(2n);\mathbb Q}(1-s)=\xi_{Sp(2n);\mathbb Q}(s);$}
\vskip 0.20cm
\noindent
(2) ({\bf The Riemann Hypothesis\,$_{Sp(2n);\mathbb Q}$})\newline
$All\ zeros\ of\ the\ zeta\ function\ \xi_{Sp(2n);\mathbb Q}(s)\ lie\ on\ the\ central\ line\
\mathrm{Re}\,s=\displaystyle{\frac{1}{2}}.$
\vskip 0.30cm
Up to this point, the importance of the period 
$$\omega_{\mathbb Q}^{G}(\lambda):=\sum_{w\in W} \Big(\frac {1}{\prod_{\alpha\in \Delta_0} \langle w\lambda-\rho, \alpha^{\vee}\rangle}  \cdot \prod_{\alpha>0, w\alpha<0} \frac {\xi\big(\langle\lambda, \alpha^{\vee}\rangle\big)}{\xi\big(\langle \lambda, \alpha^{\vee}\rangle+1\big)}\Big)$$ has been fully realized and the huge symmetry induced from the Weyl group $W$ is noticed.
\vskip 0.20cm
\noindent
{\bf A.4.2 $G_2$: Maximal Parabolics Discovered}

\smallskip
\noindent
The success of introducing natural zetas for $Sp(n)$ which are supposed to satisfying the Riemann Hypothesis proves to be very crucial. Passing this point, we then seriously try to find natural zetas for other types of classical groups. 

Practically, to be able to find such zetas, we still need solve two main technical problems:
\vskip 0.30cm
\noindent
1) how to introduce an analog of Epstein zeta function for other groups? 
Such a function should at least satisfy the property that  it can be obtained as the residue along certain singular hyperplanes of the (relative) Eisenstein series $E^{G/B}({\bf 1};\lambda)(g)$ associated to constant function one on the Borel; and
\vskip 0.30cm
\noindent
2) what are singular hyperplanes along which the residues should be taken?
\vskip 0.30cm
However, by reviewing what has been done for $SL(n)$ and $Sp(2n)$, 
for the purpose of introducing zetas, we realize that the completed theory for (1) is not really needed absolutely:  What matters (for introducing our new zetas) is not Epstein type zeta, but the period
$$\omega_{\mathbb Q}^{G}(\lambda):=\sum_{w\in W} \frac {1}{\prod_{\alpha\in \Delta_0} \langle w\lambda-\rho, \alpha^{\vee}\rangle}  \cdot \prod_{\alpha>0, w\alpha<0} \frac {\xi\big(\langle\lambda, \alpha^{\vee}\rangle\big)}{\xi\big(\langle \lambda, \alpha^{\vee}\rangle+1\big)}.$$

With (1) solved, we then shift to (2). At the very beginning, we had no idea on how to deal it -- to solve this problem we first need  understand where are singularities for $E^{G/B}({\bf 1};\lambda)(g)$; more importantly, even if knowing the singularities, we still need 
figure out along which singular hyperplanes we take the residues, as there are  many many  possible choices.

As such, at this preliminary stage of our study,  we decide to be more practical. That is, not try to solve the problem completely, but try to work with examples with the hope to expose hidden structures: 
After all, the most important points are to introduce new zetas, 
and once introduced to check whether they satisfy the functional equation and further the Riemann Hypothesis. 

For such a limited practical purpose, then clearly, among all classical groups, we need  test these groups which are with relatively smaller ranks and with 
reasonable smaller sizes of Weyl groups. By looking at $B_n,\, D_n$, $E_{6,7,8}, F_4$ and $G_2$, it is obvious why we decide to focus on $G_2$ -- $G_2$, being exceptional and interesting, is of rank two and with only 12 Weyl elements. This is extremely nice: rank two should make our study more like to be successful -- after all, the period $$\omega^{G_2}_{\mathbb Q}(z_1,z_2)=\sum_{w\in W} \frac {1}{\langle w\lambda-\rho, \alpha_{\mathrm{short}}^{\vee}\rangle\,\cdot\,\langle w\lambda-\rho, \alpha_{\mathrm{long}}^{\vee}\rangle}  \cdot \prod_{\alpha>0, w\alpha<0} \frac {\xi\big(\langle\lambda, \alpha^{\vee}\rangle\big)}{\xi\big(\langle \lambda, \alpha^{\vee}\rangle+1\big)}$$ is a  function with two variables $(z_1,z_2)=\lambda\in\frak a_0$, where $\Delta_0:=\{\alpha_{\mathrm{short}},\alpha_{\mathrm{long}}\}$ with $\alpha_{\mathrm{short}}$ the short root and $\alpha_{\mathrm{short}}$ the long root. Consequently,  we only need  find a single singular line $az_1+bz_2+c=0$. 

At this point, then by recall what has happened for $SL$ and $Sp$,
we conclude that in fact {\it all singular hyper-planes appeared there are factors
of the denominator of the term in $\omega_{\mathbb Q}^G(\lambda)$ corresponding to the identity Weyl element $\mathrm {Id}$.}
Applying this to $G_2$, we are led to $$\frac {1}{\langle w\lambda-\rho, \alpha_{\mathrm{short}}^{\vee}\rangle\,\cdot\,\langle w\lambda-\rho, \alpha_{\mathrm{long}}^{\vee}\rangle}.$$
Now it is crystal clear that we should do -- There are two possibilities for the choice of a single singular line: 

\noindent
(1) $\langle w\lambda-\rho, \alpha_{\mathrm{short}}^{\vee}\rangle=0$ or 

\noindent
(2) $\langle w\lambda-\rho, \alpha_{\mathrm{long}}^{\vee}\rangle=0$.

In this way, then we obtain two new zetas for $G_2$. Now recall that by Lie theory, there is a one-to-one and onto correspondence between maximal parabolic subgroups and simple roots, it is then only natural
for us to name the corresponding zeta functions $\xi_{\mathbb Q}^{G_2/P_{\mathrm{long}}}(s)$ and $\xi_{\mathbb Q}^{G_2/P_{\mathrm{short}}}(s)$ respectively, where $P_{\mathrm{short}}$ and $P_{\mathrm{long}}$ correspond to $\alpha_{\mathrm{long}}$ and $\alpha_{\mathrm{short}}$ respectively.
The precise calculation is carried out in Appendix B. In particular, the result confirms that we have {\it the functional equation}
$$\xi_{\mathbb Q}^{G_2/P_{\mathrm{long}}}(1-s)=\xi_{\mathbb Q}^{G_2/P_{\mathrm{long}}}(s)\qquad\mathrm{and}\qquad \xi_{\mathbb Q}^{G_2/P_{\mathrm{short}}}(1-s)=\xi_{\mathbb Q}^{G_2/P_{\mathrm{short}}}(s).$$ 
\vskip 0.30cm
\noindent
{\bf A.4.3 Zetas for $(G,P)/\mathbb Q$: Singular Hyper-planes Found}

\smallskip
\noindent
With the discovery of importance played by the period $\omega_{\mathbb Q}^G(\lambda)$ in our study of zeta functions, and the success of the discussion on $G_2$, we next want to systematically understand how singular hyperplanes are chosen in the process of taking residues. For this we go back to examine the examples of $SL(n), Sp(2n)$ and $G_2$.

\noindent
a) For $SL(n)$, a rank $(n-1)$ group, as usual, $$\Delta_0=\{e_1-e_2, e_2-e_3,\dots,e_{n-1}-e_n\},$$ with $$\lambda=(z_1,z_2,\dots,z_n)\in\frak a_0\subset \mathbb C^n, \qquad \sum_{i=1}^nz_i=0,$$ where $e_i$'s are the standard ON basis for $\mathbb C^n$. In the definition of $\xi_{SL(n), \mathbb Q}(s)$, the $(n-2)$-singular hyperplanes are chosen to be $$z_1-z_2=1, z_2-z_3=1,\dots, z_{n-2}-z_{n-1}=1;$$

\noindent
b) For $Sp(2n)$, a rank $n$ group, as usual, $$\Delta_0=\{e_1-e_2, e_2-e_3,\dots,e_{n-1}-e_n, 2e_n\}$$  with $\lambda=(z_1,z_2,\dots,z_n)\in\frak a_0=\mathbb C^n$. In the definition of $\xi_{Sp(2n), \mathbb Q}(s)$,  the $(n-1)$-singular hyperplanes are chosen to be $$z_1-z_2=1, z_2-z_3=1,\dots, z_{n-1}-z_{n}=1;$$

\noindent
c) For $G_2$, a rank two group, as usual $$\Delta_0=\{\alpha_{\mathrm{short}},\alpha_{\mathrm{long}}\}.$$ In this case, we  decided to use $\lambda=z_1(2\alpha_{\mathrm{short}}+\alpha_{\mathrm{long}})+z_2(\alpha_{\mathrm{short}}+\alpha_{\mathrm{long}})$. As said above, two different choices of a single singular line are chosen: $z_1-z_2=1$ and $z_2=0$.

As such, by looking at these singular hyperplanes more carefully, we conclude that

\noindent
a) For $SL(n)$, they are given by $$\langle w\lambda-\rho,e_1-e_2\rangle=0,\ \langle w\lambda-\rho,e_2-e_3\rangle=0, \dots, \langle w\lambda-\rho,e_{n-2}-e_{n-1}\rangle=0,$$ 
or better, are given by $$\langle w\lambda-\rho,\alpha^\vee\rangle=0,\qquad\alpha\in\Delta\backslash\{e_{n-1}-e_n\};$$
\noindent
b) For $Sp(2n)$, they are given by $$\langle w\lambda-\rho,e_1-e_2\rangle=0,\ \langle w\lambda-\rho,e_2-e_3\rangle=0, \dots, \langle w\lambda-\rho,e_{n-1}-e_n\rangle=0,$$
or better, are given by $$\langle w\lambda-\rho,\alpha^\vee\rangle=0,\qquad\alpha\in\Delta\backslash\{2e_n\};$$

\noindent
c) For $G_2$, easily with the choice $\lambda=z_1(2\alpha_{\mathrm{short}}+\alpha_{\mathrm{long}})+z_2(\alpha_{\mathrm{short}}+\alpha_{\mathrm{long}})$, the line $z_1-z_2=1$ corresponds to $\langle \lambda-\rho,\alpha_{\mathrm{short}}^\vee\rangle=0$, while
line $z_2=1$ corresponds to $\langle \lambda-\rho,\alpha_{\mathrm{long}}^\vee\rangle=0.$ Or better put, the line $z_1-z_2=1$
 is given by $$
 \langle \lambda-\rho,\alpha^\vee\rangle=0,\qquad \alpha\in\Delta_0\backslash\{\alpha_{\mathrm{long}}\};$$
while the line $z_2=1$
  is given by $$\langle \lambda-\rho,\alpha^\vee\rangle=0, \qquad\alpha\in \Delta_0\backslash\{\alpha_{\mathrm{short}}\}.$$

Recall now that, to introduce new zetas, we are determined to use
$$\omega^G_{\mathbb Q}(\lambda)=\sum_{w\in W}\Bigg(\frac {1}{\prod_{\alpha\in \Delta_0} (\langle \lambda, w^{-1}\alpha^{\vee}\rangle-1)}  \cdot \prod_{\alpha>0, w\alpha<0} \frac {\xi\big(\langle\lambda, \alpha^{\vee}\rangle\big)}{\xi\big(\langle \lambda, \alpha^{\vee}\rangle+1\big)}\Bigg),$$ a special  period governed by huge symmetries. Recall also that, for finding singular hyper-planes, our success for $SL$ and $Sp$ led to the term corresponding to $w=1$: $$\frac{1}{\prod_{\alpha\in \Delta_0}(\langle \lambda, \alpha^{\vee}\rangle-1)} \cdot 1=\frac{1}{\prod_{\alpha\in \Delta_0}(\langle \lambda, \alpha^{\vee}\rangle-1)}.$$ With such a focus, it is then not too difficult for us to detect that

\noindent
{\it all $(r-1)$-singular hyperplanes are taken from the total $r$-factors in the denominator of this term, where $r$ is the rank of the group.}

\noindent
Once this is observed, then it is extremely clear what we have done so far: a special choice of the $(r-1)$-singular hyperplanes correspond to a fixed choice of certain special maximal parabolic subgroup. More precisely,
for a fixed standard maximal parabolic subgroup $P$, by Lie theory, there exists a single simple root $\alpha_P$ such that $P$ corresponding to 
$\Delta_0\backslash\{\alpha_P\}$. As such, {\it the $(r-1)$ singular hyperplanes chosen may be understood as these given by $\langle\lambda-\rho,\alpha^\vee\rangle=0$ for $\alpha\in\Delta_0, \alpha\not=\alpha_P$.}

Upon this point, we are quite sure how a new type of zetas for $(G,P)$ should be introduced. And more importantly, we understand the importance of the role played by the symmetry. This then leads to
Definition 2 of periods of $(G,P)/\mathbb Q$:
$$\omega_{\mathbb Q}^{G/P}(s):=\\
\mathrm{Res}_{\{\langle \lambda-\rho,\alpha^\vee\rangle=0\,:\,
\alpha\in\Delta_0\backslash\{\alpha_P\}\}}\Big(\omega_{\mathbb Q}^G(\lambda)\Big)$$ where $\alpha_P$ is the simple root corresponds to the maximal parabolic $P$. With suitable normalization as done in Definition 3, we then finally obtain our new zetas $\xi_{\mathbb Q}^{G/P}(s)$ for $(G,P)$ over $\mathbb Q$, whose importance can be read from the following
\vskip 0.30cm
\noindent
{\bf Conjecture.} (1) ({\bf Functional Equation}) {\it $\xi^{G/P}_{\mathbb Q}(1-s)=\xi^{G/P}_{\mathbb Q}(s);$}
\vskip 0.20cm
\noindent
(2) ({\bf The Riemann Hypothesis\,$^{G/P}_{\mathbb Q}$})\newline
$\qquad All\ zeros\ of\ the\ zeta\ function\ \xi^{G/P}_{\mathbb Q}(s)\ lie\ on\ the\ central\ line\
\mathrm{Re}\,s=\displaystyle{\frac{1}{2}}.$

To support this new approach,
we start  working on more examples (for these new zetas) associated to other type of standard maximal subgroups (of $SL(3)$, $SL(4)$, $Sp(4)$ and $G_2$). The details are given in Appendix B.
\subsection{Conclusion Remarks}
\vskip 0.20cm
\noindent
{\bf A.5.1 Analogue of High Rank Zetas}

\smallskip
\noindent
We here propose an approach aiming at introducing  genuine zeta functions for $(G,P)/F$, as a natural generalization of high rank zeta functions. 

Denote by $\mathbb A_F$ the adelic ring of $F$.
Let $G$ be a reductive group defined over $F$, and $P$ a maximal parabolic subgroup. Then for the constant function $\bold 1$ on $P$, we form the relative Eisenstein series $E(\bold 1;\lambda_{G/P};g)=E^{G/P}(\bold 1;\lambda_{G/P};g)$. For a fixed sufficiently positive $T\in\frak a_0$, the space of characters of the Borel $B$ of $G$, introduce a single variable period
$$\omega_{G/P;F}^T(\lambda_{G/P}):=\int_{Z_{G(\mathbb A_F)}G(F)\backslash G(\mathbb A_F)}\Lambda^TE^{G/P}(\bold 1;\lambda_{G/P};g)\,d\mu(g).$$ We expect that

\noindent
An analogue of Fact G-I-J for $G$-principal lattices exists.

If so, then it makes sense to introduce $$\begin{aligned}\omega_{G/P;F}(\lambda):=&\omega_{G/P;F}^T(\lambda)|_{T=0}\\
=&\int_{\frak F_G(0)\subset Z_{G(\mathbb A_F)}G(F)\backslash G(\mathbb A_F)}E^{G/P}(\bold 1;\lambda_{G/P};g)\,d\mu(g).\end{aligned}$$ In particular, from $\omega_{G/P;F}(\lambda)$,  
a suitable normalization will then finally lead to an analogue of high rank zetas for $(G,P)/F$.

\noindent
{\bf Questions.} (1) Is it possible to get $E^{G/P}(\bold 1;\lambda_{G/P};g)$ from $E^{G/B}(\bold 1;\lambda;g)$, the relative Eisenstein series associated to the constant function $\bold 1$ on the Borel, by taking residues along with suitable $\mathrm{rank}(G)-1$ singular hyper-planes?

\noindent
(2) Can we take these singular hyper-planes simply as $\langle\lambda-\rho,\alpha^\vee\rangle=0,\ \alpha\in\Delta_0\backslash\{\alpha_P\}?$

\noindent
(3) Is it possible to introduce a completed Eisenstein series $\widehat E^{G/P}(\bold 1;\lambda_{G/P};g)$ from $E^{G/P}(\bold 1;\lambda_{G/P};g)$ so that the resulting zeta function admits only finite many singularities, satisfies a simple functional equation, and the Riemann Hypothesis?
\vskip 0.20cm
\noindent
{\bf A.5.2 $T$-version}

\smallskip
\noindent
In our discussion above, by adapting an analytic method, we can extend our discussion for periods defined originally for sufficiently positive $T$ to these for $T=0$. This makes the theory more canonical and elegant.
However the use of $T$-version proves to be quite helpful --  as example for $SL(3,4,5)$  shows, such a $T$-version can be used to help us 
to understand the additional symmetry for our new zeta functions. For example, we know that $$\xi_{\mathbb Q}^{SL(3)/P_{2,1}}(s)=\xi_{\mathbb Q}^{SL(3)/P_{1,2}}(s),\qquad\xi_{\mathbb Q}^{SL(4)/P_{3,1}}(s)=\xi_{\mathbb Q}^{SL(4)/P_{1,3}}(s),$$ and $$\xi_{\mathbb Q}^{SL(5)/P_{4,1}}(s)=\xi_{\mathbb Q}^{SL(5)/P_{1,4}}(s),\qquad\xi_{\mathbb Q}^{SL(5)/P_{2,3}}(s)=\xi_{\mathbb Q}^{SL(5)/P_{3,2}}(s).$$ On surface, these relations may be viewed as a reflection of the symmetry between the Eisenstein series $E_{r-m,m}$ associated to the maximal parabolic $P_{r-m,m}$ and the Eisenstein series $E_{m, r-m}$ associated to the maximal parabolic $P_{m,r-m}$. (See A.3.1 for details.)
More deeply, it roots into the symmetry between  $P_{r-m,m}$ and $P_{m,r-m}$ for maximal parabolic subgroups of $SL(r)$. 

Put this in concrete term, for $SL(3)$, we can further introduce  $T$-version 
zeta functions $\xi_{\mathbb Q}^{SL(3)/P_{2,1};T}(s)$  and $\xi_{\mathbb Q}^{SL(3)/P_{1,2};T}(s)$, analogues of  $\xi_{\mathbb Q}^{SL(3)/P_{2,1}}(s)$  and $\xi_{\mathbb Q}^{SL(3)/P_{1,2}}(s)$ respectively, starting from the $T$-version period $\omega_{\mathbb Q}^{SL(3);T}(\lambda)$ in A.2.6. Then one checks that with $T\in\mathbb C\cdot\rho$, i.e., with $T$ specialized as points on the line spanned by $\rho$, we have 
$$\xi_{\mathbb Q}^{SL(3)/P_{2,1};T}(1-s)=\xi_{\mathbb Q}^{SL(3)/P_{1,2};T}(s).$$ This is then the root of the equality
$$\xi_{\mathbb Q}^{SL(3)/P_{2,1}}(s)=\xi_{\mathbb Q}^{SL(3)/P_{1,2}}(s).$$ We expect that holds for all zetas related to $(SL(r),P_{r-m,m})/\mathbb Q.$

Along with this line, then we also expect that the symmetry, or better, the duality, between type $B_n$ and $C_n$ groups will have similar impact to our new zetas. In a sense, various symmetries are the main reason why our new zetas satisfy the functional equations and the Riemann Hypothesis.

We end this $T$-version discussion by pointing out that the Riemann Hypothesis does not hold for $\xi_{\mathbb Q}^{SL(3)/P_{2,1};T}(s)$ if $T$ is not 0. So our new zetas $\xi_{\mathbb Q}^{G/P}(s)$, being specialization of $T$-version zetas $\xi_{\mathbb Q}^{G/P;T}(s)$ to the ground zero and hence delicate, 
are quite canonical, hence absolutely beautiful.

\vskip 0.20cm
\noindent
{\bf A.5.3 Where Lead To}

\smallskip
\noindent
It is hard to predict, being new and rich. In general terms, two aspects are worth being mentioned. One is for the zetas themselves, the other is for possible applications.

For zetas themselves, the first and the up-most task is then concentrated on the (proof of) functional equations and the corresponding Riemann Hypothesis. Examples listed in Appendix B for $SL(2,3,4,5), Sp(4)$ and $G_2$ show that the associated zetas satisfy the Functional Equation. This is beautiful,  reflecting additional symmetry, and supposedly doable even expected to be very complicated. On the other hand, for the  Riemann Hypothesis associated to new zetas, responding to our inquires  ([W4]), Suzuki first made several crucial numerical tests on zeros of zetas $\xi_{SL(4);\mathbb Q}(s)$, $\xi_{SL(5);\mathbb Q}(s)$ and $\xi_{Sp(4);\mathbb Q}(s)$ ([S2]). Shortly after, in January 2008, he was able to theoretically verify the Riemann Hypothesis for zetas $\xi_{Sp(4);\mathbb Q}(s)$ and  $\xi_{\mathbb Q}^{G_2/P}(s)$ ([S3,\,4]), by strengthening a method used for establishing the RH of $\xi_{SL(2);\mathbb Q}(s)$ ([LS])
and of $\xi_{SL(3);\mathbb Q}(s)$ ([S]). (In fact, this method can also be used to show that outside a certain finite box, all zeros of $\xi_{\mathbb Q}^{Sp(4)/P_{2e_2}}(s)$ lie on the line $\mathrm{Re}\,(s)=\frac{1}{2}$ as well.)

The third is about a generalization to all reductive groups. Even physically, this can be done simply since all the framework works in this generality. But we are somehow a bit hesitated feeling
that time is not ripe to make such a move, even we know that, up to a constant factor,\newline
\centerline{$\xi_F^{G_1\times G_2/P_1\times G_2}(s)=\xi_F^{G_1/P_1}(s)$} \newline
and that the RH holds for all rank 2 groups (modulo the finite box mentioned above for $\xi_{\mathbb Q}^{Sp(4)/P_{2e_2}}(s)$).
\vskip 0.30cm
For applications,  an obvious is about the relation between new zetas and the classical Riemann zeta function. Problems likely to be asked here are: what should be the relations between their zeros?  This can be put more precisely, for example, as: if we just consider a series, e.g., the series for $SL(r)/P_{r-1,1}$, or a collection, e.g., the collection of rank $r$ groups, what should be the sequence of the $n$-th zeros for a fixed $n$? what about the distributions of these zeros, the gaps between ordered pairs of zeros? etc. For this, a related interesting point should be mentioned: the completed Riemann zeta function can be written as a difference between  two entire functions which both satisfy the RH. This is a new structure emerged in our understanding of $\xi_{SL(3);\mathbb Q}(s)$. (See also [S3] for  $\xi_{Sp(4);\mathbb Q}(s)$.)

We end this appendix by proposing a bit indirect, but quite speculating use
of our new zetas. We call this a \lq wonderful idea' -- the final gold is to replace the original Riemann Hypothesis in the study of distribution of primes, of classical problems such as the Goldbach conjecture, etc., with the RH for our zetas, some of which have been established.
\eject
\vskip 1.0cm
\vskip 0.45cm
\centerline {\bf\Large{REFERENCES}}
\vskip 0.20cm 
\noindent
[Ar1] J. Arthur, A trace formula for reductive groups. I. 
Terms associated to classes in $G({\mathbb Q})$. Duke Math. J. {\bf 45} (1978), 
no. 4, 911--952
\vskip 0.20cm
\noindent
[Ar2] J. Arthur, A trace formula for reductive groups. II. Applications of a 
truncation operator. Compositio Math. {\bf 40} (1980), no. 1, 87--121.
\vskip 0.20cm
\noindent
[Ar3] J. Arthur, A measure on the unipotent variety, Canad. J. Math {\bf 37},
(1985) pp. 1237--1274
\vskip 0.20cm
\noindent
[Bo1] A. Borel, Some finiteness properties of adele groups over
number fields, Publ. Math., IHES, {\bf 16} (1963) 5-30
\vskip 0.20cm
\noindent
[Bo2] A. Borel, {\it Introduction aux groupes arithmetictiques}, Hermann, 1969
\vskip 0.20cm 
\noindent  
{[Co]} A. Connes, Trace formula in noncommutative geometry and the zeros 
of the Riemann zeta function.  Selecta Math. (N.S.) 5  (1999),  no. 1, 29--106.
\vskip 0.20cm
\noindent
[D] B. Diehl, Die analytische Fortsetzung der Eisensteinreihe zur Siegelschen Modulgruppe, J. reine angew. Math., 317 (1980) 40-73
\vskip 0.20cm 
\noindent  
{[De1]} C. Deninger,  On the $\Gamma$-factors attached to motives.
Invent. Math. {\bf 104} (1991), no. 2, 245--261.
\vskip 0.20cm 
\noindent  
{[De2]} C. Deninger,  Local $L$-factors of
motives and regularized determinants.  Invent. Math.  {\bf 107}  (1992),  
no. 1, 135--150.
\vskip 0.20cm 
\noindent  
{[De3]} C. Deninger, Motivic $L$-functions and regularized determinants, 
Motives (Seattle, WA, 1991), 707-743, Proc. Sympos. Pure Math, {\bf 55} Part 1,
AMS, Providence, RI, 1994
\vskip 0.20cm 
\noindent  
{[E]} H.M. Edwards, {\it Riemann's Zeta Function}, Dover Pub., 1974
\vskip 0.20cm
\noindent 
[GS] G. van der Geer \& R. Schoof, Effectivity of Arakelov
Divisors and the Theta Divisor of a Number Field, Sel. Math., New ser.
{\bf 6} (2000), 377-398 
\vskip 0.20cm 
\noindent  
{[Gr1]} D.R. Grayson, Reduction theory using semistability. Comment. 
Math. Helv. {\bf 59} (1984), no. 4, 600--634
\vskip 0.20cm 
\noindent  
{[Gr2]} D.R. Grayson, Reduction theory using semistability. II.  
Comment. Math. Helv. {\bf 61}  (1986),  no. 4, 661--676. 
\vskip 0.20cm
\noindent
[Ha] T. Hayashi, Computation of Weng's rank 2 zeta function over an algebraic number field, J. Number Theory {\bf 125} (2007), no. 2, 473--527
\vskip 0.20cm 
\noindent
[Hu] J. Humphreys, {\it Introduction to Lie algebras and representations}, Springer-Verlag, 1972
\vskip 0.20cm 
\noindent  
{[Iw]}  K. Iwasawa, Letter to Dieudonn\'e, April 8, 1952, in
{\it Zeta Functions in Geometry}, Advanced
Studies in Pure Math. {\bf 21} (1992), 445-450
\vskip 0.20cm
\noindent 
[JLR] H. Jacquet, E. Lapid \& J. Rogawski,  Periods of automorphic forms. 
J. Amer. Math. Soc. 12 (1999), no. 1, 173--240
\vskip 0.20cm 
\noindent  
{[Ki]} H. Ki, All but finitely many non-trivial zeros of the approximations of 
the Epstein zeta function are simple and on the critical line. Proc. London Math. 
Soc. (3) 90 (2005), no. 2, 321--344. 
\vskip 0.20cm 
\noindent  
[KW] H.H. Kim \& L. Weng, Volume of truncated fundamental domains, Proc. AMS, Vol. 135 (2007) 1681-1688
\vskip 0.20cm
\noindent
[Laf] L. Lafforgue, {\it Chtoucas de Drinfeld et conjecture de 
Ramanujan-Petersson}. Asterisque No. 243 (1997)
\vskip 0.20cm
\noindent
[LS] J. Lagarias \& M. Suzuki, The Riemann Hypothesis for certain 
integrals of Eisenstein series, Journal of Number Theory, 118(2006) 98-122
\vskip 0.20cm 
\noindent  
{[L1]} S. Lang, {\it Algebraic Number Theory}, 
Springer-Verlag, 1986
\vskip 0.20cm 
\noindent  
{[L2]}  S. Lang, {\it Introduction to Arakelov theory}, Springer Verlag, 1988
\vskip 0.20cm
\noindent
[La1] R. Langlands, {\it On the functional equations satisfied by 
Eisenstein series}, Springer LNM {\bf 544}, 1976
\vskip 0.20cm
\noindent
[La2] R. Langlands, the volume of the fundamental domain for some arithmetical 
subgroups of Chevalley groups, in {\it Algebraic Groups and Discontinuous 
Subgroups,} Proc. Sympos. Pure Math. 9, AMS (1966) pp.143--148
\vskip 0.20cm 
\noindent  
[La3] R. Langlands,  {\it Euler products}, Yale Math. Monograph, Yale Univ. Press, 1971
\vskip 0.20cm 
\noindent  
{[Mi]} H. Minkowski, {\it Geometrie der Zahlen}, Leipzig and Berlin, 1896
\vskip 0.20cm
\noindent
[MW] C. Moeglin \& J.-L. Waldspurger, {\it Spectral decomposition 
and Eisenstein series}. Cambridge Tracts in Math, {\bf 113}. 
Cambridge University Press, 1995
\vskip 0.20cm 
\noindent  
{[M]} D. Mumford, {\it Geometric Invariant Theory}, Springer-Verlag, 
 (1965)
\vskip 0.20cm 
\noindent    
{[NS]} M.S.  Narasimhan \& C.S.  Seshadri,  Stable 
and unitary vector bundles on a compact Riemann surface. Ann. of Math. (2) 
82 1965
\vskip 0.20cm 
\noindent  
{[Ne]} J. Neukirch, {\it Algebraic Number Theory}, Grundlehren der
Math. Wissenschaften, Vol. {\bf 322}, Springer-Verlag, 1999
\vskip 0.20cm 
\noindent  
{[RR]} S. Ramanan \& A. Ramanathan,  Some remarks on the instability 
flag.  Tohoku Math. J. (2)  36  (1984),  no. 2, 269--291.
\vskip 0.20cm 
\noindent  
{[Ser]} J.-P. Serre, {\it Algebraic Groups and Class Fields}, GTM 117, 
Springer (1988)
\vskip 0.20cm
\noindent 
[Sie] C.L. Siegel, {\it Lectures on the geometry of numbers}, notes by B.
Friedman, rewritten by K. Chandrasekharan with the assistance of R. Suter, 
Springer-Verlag, 1989.
\vskip 0.20cm
\noindent   
{[St1]} U. Stuhler,   Eine Bemerkung zur Reduktionstheorie quadratischer 
Formen, Arch. Math. (Basel) {\bf 27} (1976), no. 6, 604--610
\vskip 0.20cm
\noindent   
{[St2]} U. Stuhler, Zur Reduktionstheorie der positiven quadratischen Formen. 
II, Arch. Math. (Basel)  {\bf 28}  (1977), no. 6, 611--619
\vskip 0.20cm
\noindent
[S] M. Suzuki, A proof of the Riemann Hypothesis for the Weng zeta function of rank 3 for the rationals, pp.175-200, in {\it Conference on L-Functions}, World Sci. 2007
\vskip 0.20cm
\noindent
[S2] M. Suzuki, private communications, Oct.-Dec., 2007
\vskip 0.32cm
\noindent
[S3] M. Suzuki, The Riemann hypothesis for Weng's zeta function of ${\rm Sp}(4)$ over $\mathbb Q$, with an appendix [W5], preprint, 2008
\vskip 0.32cm
\noindent
[SW] M. Suzuki \& L. Weng, Zeta functions for $G_2$ and their zeros,  preprint, 2008
\vskip 0.20cm 
\noindent  
{[T]} J. Tate, Fourier analysis in number fields and Hecke's
zeta functions, Thesis, Princeton University, 1950 
\vskip 0.20cm 
\noindent  {[Te]} A. Terras, {\it Harmonic analysis on symmetric spaces and applications II}, Springer-Verlag, 1988
\vskip 0.20cm 
\noindent  {[Ve]} A.B. Venkov, On the trace formula for $SL(3,\mathbb Z)$, J Soviet Math., 12 (1979), 384-424
\vskip 0.20cm 
\noindent  
{[We]}  A. Weil, {\it Basic Number Theory}, Springer-Verlag, 1973
\vskip 0.20cm
\noindent
[W-3] L. Weng, Analytic truncation and Rankin-Selberg versus algebraic
truncation and non-abelian zeta, Algebraic Number Theory and Related Topics, RIMS Kokyuroku, No.1324 (2003), 7-21.
\vskip 0.20cm
\noindent 
{[W-2]} L. Weng, Rank Two Non-Abelian Zeta and its Zeros, \newline
 available at http://xxx.lanl.gov/abs/math.NT/0412009
\vskip 0.20cm
\noindent  
{[W-1]} L. Weng, Automorphic Forms, Eisenstein Series and Spectral 
Decompositions, 
{\it Arithmetic Geometry and Number Theory}, 123-210, World Sci. 2006
\vskip 0.20cm
\noindent
[W0] L. Weng, Non-abelian zeta function for function fields, Amer. J. Math 127 (2005), 973-1017
\vskip 0.20cm
\noindent
[W1] L. Weng, Geometric Arithmetic: A Program, in {\it Arithmetic Geometry and Number Theory}, pp. 211-390, World Sci. (2006)
\vskip 0.20cm
\noindent
[W2] L. Weng,  A Rank two zeta and its zeros, J of Ramanujan Math. Soc, 21 (2006), 205-266
\vskip 0.30cm
\noindent
[W3] L. Weng, A geometric approach to $L$-functions, in {\it Conference on L-Functions}, pp. 219-370, World Sci (2007)
\vskip 0.20cm
\noindent
[W4] L. Weng, Zetas for $SL(4), SL(5), Sp(4),$ and $G_2$ over $\mathbb Q$, private notes to  Suzuki, Oct.-Dec. 2007
\vskip 0.20cm
\noindent  
{[W5]} L. Weng, Zeta function for $Sp(2n)$, 
Appendix to [S3],  preprint, 2008
\vskip 0.20cm
\noindent
[Z] Zagier, D.  The Rankin-Selberg method for automorphic functions which are not of rapid decay. 
J. Fac. Sci. Univ. Tokyo Sect. IA Math. 28(3), 415--437 (1982)
\eject
\section{\bf\Large Examples}
\vskip 0.30cm
We here list  zetas $\xi_{\mathbb Q}^{G/P}$  for $G=\,SL(2,\,3,\,4,\,5)$, $Sp(4)$ and $G_2$. Consequence,  all these zetas satisfy the  FE $\ \xi^{G/P}_{\mathbb Q}(1-s)=\xi^{G/P}_{\mathbb Q}(s).$
(Detailed calculations were given in version 2007 of this paper, but are omitted here as zetas for $SL(2,\,3)$, $Sp(4)$ and $G_2$ are now available in [W1,\,3,\,4] and [SW] respectively).

\bigskip
\noindent
{\bf Contents
\vskip 0.20cm
\noindent
\begin{footnotesize}
B.1 $SL(n)$
\vskip 0.15cm
B.1.1 $SL(2)$
\vskip 0.15cm
B.1.2 $SL(3)$
\vskip 0.15cm
B.1.3 $SL(4)$
\vskip 0.15cm
B.1.4 $SL(5)$
\vskip 0.15cm
\noindent
B.2 $Sp(4)$
\vskip 0.15cm
\noindent
B.3 $G_2$
\vskip 0.15cm
\noindent
B.4 $T$-version for $SL(3)$
\end{footnotesize}}

\subsection{$SL(n)$}
\subsubsection{$SL(2)$}
A degenerate case, since  $P=B$, the Borel. We have 
 \begin{equation}
\boxed{\xi^{SL(2)/B}_{\mathbb Q}(s)=\xi_{\mathbb Q,2}:(s)=\frac{\xi_{\mathbb Q}(2s)}{s-1}-
\frac{\xi_{\mathbb Q}(2s-1)}{s}}
\end{equation}
It is the first natural example exposed that satisfies the RH ([W1,2,3], [LS]).

\subsubsection{$SL(3)$}
Two maximal parabolic subgroups $P$, corresponding to partitions $3=2+1=1+2$. They share the same zetas:
\begin{equation}
\boxed{\begin{aligned}\xi^{SL(3)/P}_{\mathbb Q}(s)=&\xi_{\mathbb Q}(2)\cdot\frac{1}{3s-3}\cdot\xi_{\mathbb Q}(3s)\\
&-\xi_{\mathbb Q}(2)\cdot\frac{1}{3s}\cdot\xi_{\mathbb Q}(3s-2)\\
&+\frac{1}{3}\cdot\frac{1}{3s-3}\cdot\xi_{\mathbb Q}(3s-1)\\
&-\frac{1}{3}\cdot\frac{1}{3s}\cdot\xi_{\mathbb Q}(3s-1)\\
&+\frac{1}{2}\cdot\frac{1}{3s-1}\cdot\xi_{\mathbb Q}(3s-2)\\
&-\frac{1}{2}\cdot\frac{1}{3s-2}\cdot\xi_{\mathbb Q}(3s)\end{aligned}}
\end{equation} 
Contradicting to the claim in Ch.\,9 of [W3], by examining poles,
this example shows that $\xi_{\mathbb Q,r}(s)\not=\xi_{\mathbb Q}^{SL(r)/P_{r-1,1}}(s).$ So
{\it high rank zetas $\xi_{F,r}(s)$ are different from the zetas for $(SL(r),P_{r-1,1})/\mathbb Q$}. The RH is confirmed ([S]).

\subsubsection{$SL(4)$}
Three maximal parabolic subgroups $P$, corresponding to partitions $4=3+1=2+2=1+3$. Denote the corresponding maximal parabolic subgroups by $P_{3,1}, P_{2,2}, P_{1,3}$ respectively. We now know that $P_{1,3}$ and $P_{3,1}$ share the same zetas, while the zeta for $P_{2,2}$ is different. More precisely, they read  as follows:
\begin{equation}
\boxed{\begin{aligned}
\xi^{SL(4)/P_{3,1}}_{\mathbb Q}(s)=&\xi^{SL(4)/P_{1,3}}_{\mathbb Q}(s)\\
=&\frac{1}{4s-4}\xi(2)\xi(3)\cdot\xi(4s)-\frac{1}{4s}\xi(2)\xi(3)\cdot\xi(4s-3)\\
&+\frac{1}{4}\frac{1}{4s-2}
\cdot\xi(4s)-\frac{1}{4}\frac{1}{4s-2}\cdot\xi(4s-3)\\
&+\frac{1}{3}\Big[\frac{1}{4s-1}+\frac{1}{4s-2}\Big]\xi(2)\cdot \xi(4s-3)\\
&-\frac{1}{3}\Big[\frac{1}{4s-2}+\frac{1}{4s-3}\Big]\xi(2)\cdot\xi(4s)\\
&+\frac{1}{2}\frac{1}{(4s)(4s-3)}\cdot\xi(4s-1)\\
&
+\frac{1}{2}\frac{1}{(4s-1)(4s-4)}\cdot\xi(4s-2)\\
&-\frac{1}{(4s)(4s-4)}\xi(2)\cdot\xi(4s-1)\\
&-\frac{1}{(4s)(4s-4)}\xi(2)\cdot\xi(4s-2)\end{aligned}}
\end{equation}
and
\begin{equation}
\boxed{\begin{aligned}
\xi^{SL(4)/P_{2,2}}_{\mathbb Q}(s)
&:=\frac{1}{2s-3}\xi(2)^2\cdot\xi(2s)\xi(2s+1)-\frac{1}{2s+1}\xi(2)\cdot\xi(2s-2)\xi(2s-1)\\
&+\frac{1}{2s-1}\cdot\frac{1}{4}\cdot\xi(2s)\xi(2s+1)-\frac{1}{2s-1}\cdot\frac{1}{4}\cdot\xi(2s-2)\xi(2s-1)\\
&+\frac{1}{(2s)^2(2s-3)}\cdot\xi(2s-1)^2-\frac{1}{(2s-2)^2(2s+1)}\cdot\xi(2s)^2\\
&-\frac{1}{2s-2}\xi(2)\cdot\xi(2s)\xi(2s+1)+\frac{1}{2s}\xi(2)\cdot\xi(2s-2)\xi(2s-1)
\\
&+\frac{1}{(2s-2)(2s)}\cdot\xi(2s-1)\xi(2s)\\
&-
\frac{2}{(2s-3)(2s+1)}\xi(2)\cdot\xi(2s-1)\xi(2s)\end{aligned}}
\end{equation}

\subsubsection{$SL(5)$}
Four maximal parabolic subgroups correspond to  the partitions $5=4+1=3+2=2+3=1+4$. Denote the associated standard maximal parabolic subgroups by $P_{4,1}, P_{3,2}, P_{2,3}, P_{1,4}$ respectively.
Then we know that the zeta for $P_{4,1}$ is the same as that for $P_{1,4}$, while 
the zeta for $P_{3,2}$ is the same as that for $P_{2,3}$.

More precisely, the new zeta functions $\xi_{\mathbb Q}^{SL(5)/P_{4,1}}(s)=\xi_{\mathbb Q}^{SL(5)/P_{1,4}}(s)$ are given by
 \begin{equation}
 \boxed{\begin{aligned}&\xi_{\mathbb Q}^{SL(5)/P_{4,1}}(s)=\xi_{\mathbb Q}^{SL(5)/P_{1,4}}(s)=\xi_{SL(5);\mathbb Q}(s):=\\
&\Big[\frac{1}{5s-5}\xi(5s)-\frac{1}{5s}\xi(5s-4)\Big]\xi(2)\xi(3)\xi(4)\\
&+\frac{1}{4}\Big\{\Big[\frac{1}{5s-1}\xi(5s-4)-\frac{1}{5s-4}\xi(5s)\Big]\\
&+\Big[\frac{1}{5s-3}\xi(5s-4)-\frac{1}{5s-2}\xi(5s)\Big]\Big\}\xi(2)\xi(3)\\
&+\frac{1}{9}\Big[\frac{1}{5s-2}\xi(5s)-\frac{1}{5s-3}\xi(5s-4)\Big]\xi(2)\\
&+\frac{1}{6}\Big\{\Big[\frac{1}{5s-3}\xi(5s)-\frac{1}{5s-2}\xi(5s-4)\Big]\\
&+\Big[\frac{1}{5s-2}\xi(5s)-\frac{1}{5s-3}\xi(5s-4)\Big]\Big\}\xi(2)\\
&+\Big\{\frac{1}{3}\Big[\frac{1}{5s(5s-4)}\xi(5s-1)+\frac{1}{(5s-5)(5s-1)}\xi(5s-3)\Big]\\
&+\frac{1}{2}\Big[\frac{1}{(5s-1)(5s-5)}\xi(5s-2)+\frac{1}{(5s-4)(5s)}\xi(5s-2)\Big]\\
&+\frac{1}{3}\Big[\frac{1}{(5s-2)(5s-5)}\xi(5s-3)+\frac{1}{(5s-3)(5s)}\xi(5s-1)\Big]\Big\}\xi(2)\\
&+\frac{1}{8}\Big[\frac{1}{5s-3}\xi(5s-4)-\frac{1}{5s-2}\xi(5s)\Big]\\
&+\frac{1}{4}\Big[\frac{1}{5s-2}\xi(5s-4)-\frac{1}{5s-3}\xi(5s)\Big]\xi(2)^2\\
&-\frac{1}{4}\Big[\frac{1}{(5s-3)(5s)}\xi(5s-1)+\frac{1}{(5s-2)(5s-5)}\xi(5s-3)\Big]\\
&-\Big[\frac{1}{(5s)(5s-5)}\xi(5s-1)+\frac{1}{(5s)(5s-5)}\xi(5s-3)\Big]\xi(2)\xi(3)\\
&-\frac{1}{4}\frac{1}{(5s-1)(5s-4)}\xi(5s-2)-\frac{1}{(5s)(5s-5)}\xi(5s-2)\xi(2)^2\end{aligned}}
\end{equation}
(which, as well as the next one, is quite complicated to obtain: totally 1200 cases should be discussed from which further residues should be taken,)
and
\begin{equation}
\boxed{\begin{aligned}
&\xi^{SL(5)/P_{3,2}}_{\mathbb Q}(s+1):=\xi^{SL(5)/P_{3,2}}_{\mathbb Q;o}(s)=\xi^{SL(5)/P_{2,3}}_{\mathbb Q;o}(s)\\
=:&\frac{1}{5s}\xi(2)^2\xi(3)\cdot\xi(5s+4)\xi(5s+5)+\frac{1}{4(5s+2)}\xi(2)\cdot\xi(5s+4)\xi(5s+5)\\
&+\frac{1}{(5s+4)^2(5s)}\xi(2)\cdot \xi(5s+2)\xi(5s+3)-\frac{1}{2(5s+1)}\xi(2)\xi(3)\cdot\xi(5s+4) \xi(5s+5)\\
&-\frac{1}{3(5s+1)}\xi(2)^2\cdot\xi(5s+4)\xi(5s+5)-\frac{1}{3(5s+2)}\xi(2)^2\cdot\xi(5s+4)\xi(5s+5)\\
-&\frac{1}{4(5s+2)(5s+4)}\cdot\xi(5s+3)\xi(5s+4)+\frac{1}{3(5s+3)}\xi(2)^2\cdot\xi(5s+1)\xi(5s+2)\\
&-\frac{1}{2(5s+1)(5s+3)(5s+4)}\cdot \xi(5s+2)\xi(5s+3)+\frac{1}{8(5s+2)}\cdot\xi(5s+1) \xi(5s+2)\\
&-\frac{1}{(5s+1)^2(5s+5)}\xi(2)\cdot\xi(5s+3)\xi(5s+4)-\frac{1}{6(5s+3)}\xi(2)\cdot\xi(5s+1)\xi(5s+2)\\
&-\frac{1}{2(5s)(5s+3)^2}\cdot\xi(5s+2)^2-\frac{1}{4(5s+2)(5s+3)}\cdot\xi(5s+2)\xi(5s+4)\\
&+\frac{1}{2(5s+4)}\xi(2)\xi(3)\cdot \xi(5s+1)\xi(5s+2)+\frac{1}{2(5s)(5s+4)}\xi(2)\cdot\xi(5s+2) \xi(5s+3)\\
&+\frac{1}{2(5s+1)(5s+4)}\xi(2)\cdot\xi(5s+2)\xi(5s+3)+\frac{1}{6(5s+2)}\xi(2)\cdot\xi(5s+4)\xi(5s+5)\\
&+\frac{1}{6(5s+3)}\xi(2)\cdot\xi(5s+4)\xi(5s+5)+\frac{1}{2(5s+1)(5s+4)}\xi(2)\cdot\xi(5s+3)\xi(5s+4)\\
&+\frac{1}{(5s+1)^2(5s+4)^2}\cdot \xi(5s+3)^2+\frac{1}{3(5s+2)(5s+4)}\xi(2)\cdot\xi(5s+2) \xi(5s+4)\\
&-\frac{1}{8(5s+3)}\cdot\xi(5s+4)\xi(5s+5)-\frac{1}{6(5s+2)}\xi(2)\cdot\xi(5s+1)\xi(5s+2)\\
&-\frac{1}{4(5s+3)}\xi(2)\cdot\xi(5s+1)\xi(5s+2)+\frac{1}{2(5s+1)(5s+5)}\xi(2)\cdot\xi(5s+3)\xi(5s+4)\\
&+\frac{1}{2(5s+2)^2(5s+5)}\cdot \xi(5s+4)^2-\frac{1}{(5s)(5s+5)}\xi(2)\xi(3)\cdot\xi(5s+2) \xi(5s+4)\\
&-\frac{1}{(5s)(5s+5)}\xi(2)^2\cdot\xi(5s+3)\xi(5s+4)-\frac{1}{(5s+1)(5s+2)(5s+5)}\xi(2)\cdot\xi(5s+4)^2\\
&+\frac{1}{3(5s+1)(5s+3)}\xi(2)\cdot\xi(5s+2)\xi(5s+4)\\
&+\frac{1}{2(5s+1)(5s+2)(5s+4)}\cdot\xi(5s+3)\xi(5s+4)\\
&-\frac{1}{4(5s+1)(5s+3)}\cdot \xi(5s+2)\xi(5s+3)-\frac{1}{(5s+5)}\xi(2)^2\xi(3)\cdot\xi(5s+1) \xi(5s+2)\\
&-\frac{1}{(5s)(5s+5)}\xi(2)^2\cdot\xi(5s+2)\xi(5s+3)+\frac{1}{(5s)(5s+3)(5s+4)}\xi(2)\cdot\xi(5s+2)^2\\
&+\frac{1}{3(5s+4)}\xi(2)^2\cdot\xi(5s+1)\xi(5s+2)\end{aligned}}
\end{equation}

\subsection{$Sp(4)$}

Two maximal parabolic subgroups corresponding to simple roots $\{e_1-e_2\}$ and $\{2e_2\}$  respectively. Their zetas read as follows:
\begin{equation}
\boxed{\begin{aligned}
\xi^{Sp(4)/P_{e_1-e_2}}_{\mathbb Q}(s)=&\frac{1}{s-2}\xi(2)\cdot\xi(s+1)\xi(2s)-\frac{1}{s+1}\xi(2)\cdot\xi(s-1)\xi(2s-1)\\
&-\frac{1}{2s-2}\cdot \xi(s+1)\xi(2s)+\frac{1}{2s}\cdot\xi(s-1) \xi(2s-1)\\
&-\frac{1}{(2s-2)(s+1)}\cdot\xi(s)\xi(2s)-\frac{1}{(2s)(s-2)}\cdot\xi(s)\xi(2s-1)\end{aligned}}
\end{equation}
and
\begin{equation}
\boxed{\begin{aligned}
\xi^{Sp(4)/P_{2e_2}}_{\mathbb Q}(s)=&\frac{1}{2s-3}\xi(2)\cdot\xi(2s+1)-\frac{1}{2s+1}\xi(2)\cdot\xi(2s-2)\\
&-\frac{1}{2(2s-1)}\cdot\xi(2s+1)+\frac{1}{2(2s-1)}\cdot\xi(2s-2)\\
&-\frac{1}{(2s+1)(2s-2)}\cdot\xi(2s)-\frac{1}{(2s)(2s-3)}\cdot\xi(2s-1).\end{aligned}}
\end{equation}
The RH for $\xi_{\mathbb Q}^{Sp(4)/P_1}(s)$ is confirmed  ([S2]), whose method, a generalization of ([S] and/or [SW]), can also be used to show that
outside a finite box, all zeros of $\xi_{\mathbb Q}^{Sp(4)/P_2}(s)$ lie on the line $\mathrm{Re}(s)=\frac{1}{2}.$
\subsection{$G_2$}
Two maximal parabolic subgroups corresponding to the long and the short root respectively. Their zetas read as follows:
 \begin{equation}
\boxed{\begin{aligned}
\xi^{G_2/P_{\mathrm{long}}}_{\mathbb Q}(s)=&\frac{1}{s-2}\xi(2)\cdot\xi(s+1)\xi(2s)\xi(3s)\\
&-\frac{1}{s+1}\xi(2)\cdot\xi(s-1)\xi(2s-1)\xi(3s-2)\\
&-\frac{1}{2s-2}\cdot\xi(s+1)\xi(2s)\xi(3s)+\frac{1}{2s}\cdot\xi(s-1)\xi(2s-1)\xi(3s-2)\\
&-\frac{1}{(3s)(2s-2)}\cdot\xi(s)\xi(2s)\xi(3s-1)\\
&-\frac{1}{(3s-1)(s-2)}\cdot\xi(s)\xi(2s-1)\xi(3s-2)\\
&-\frac{1}{(3s-3)(2s)}\cdot\xi(s)\xi(2s-1)\xi(3s-1)\\
&-\frac{1}{(3s-2)(s+1)}\cdot\xi(s)\xi(2s)\xi(3s)\end{aligned}}
\end{equation}
and
\begin{equation}
\boxed{\begin{aligned}
\xi^{G_2/P_{\mathrm{short}}}_{\mathbb Q}(s)=&\frac{1}{s-3}\xi(2)\cdot\xi(s+2)\xi(2s)-\frac{1}{s+2}\xi(2)\cdot\xi(s-2)\xi(2s-1)\\
&+\frac{1}{2s-2}\cdot\xi(s-2)\xi(2s-1)-\frac{1}{2s}\cdot\xi(s+2)\xi(2s)\\
&-\frac{1}{s(s-3)}\cdot\xi(s-1)\xi(2s-1)-\frac{1}{(s-1)(s+2)}\cdot\xi(s+1)\xi(2s)\\
&-\frac{1}{(2s-2)(s+1)}\cdot\xi(s)\xi(2s)
-\frac{1}{(2s)(s-2)}\cdot\xi(s)\xi(2s-1)\end{aligned}}
\end{equation}

The RH for $\xi_{\mathbb Q}^{G_2/P}(s)$ is confirmed by Suzuki ([SW]).
\subsection{$T$-Version for $SL(3)$}

In this subsection, we indicate how functional equation for our zetas can be obtained from a general $T$-construction. For simplicity, we consider only $G=SL(3)$.

By definition, $$\omega_{\mathbb Q}^{G;T}(s)=\sum_{w\in W}\Bigg(\frac{\langle w\lambda-\rho,T\rangle}{\prod_{\alpha\in\Delta_0}
\langle w\lambda-\rho,\alpha^\vee\rangle}\cdot\prod_{\alpha>0,w\alpha<0}\frac{\xi(\langle \lambda,\alpha^\vee\rangle)}
{\xi(\langle \lambda,\alpha^\vee\rangle+1)}\Bigg).$$ In particular, for $G=SL(3)$, we may take $\lambda=(z_1,z_2,z_3)$ with $z_1+z_2+z_3=0$, $T=(x,y,-x-y)$, $\rho=(1,0,-1)$ and $W=S_3$ with $w\in W=S_3$ acts via the corresponding permutation on lower indices.

Thus by taking residue along $z_1-z_2=1$ and assuming $z_2=t, z_1=t+1, z_3=-2t-1$, we get, using the tables in subsection B.1.2,
$$\begin{aligned}\omega_{\mathbb Q}^{SL(3)/P_{1,2};T}(t)=&\frac{1}{3t}\xi(2)\cdot \xi(3t+3)\cdot e^{3tx+3ty+4x+2y}\\
&-\frac{1}{2}\frac{1}{3t+1}\cdot\xi(3t+3)\cdot e^{(3t+3)(x+y)}\\
&+\frac{1}{2}\frac{1}{3t+2}\cdot\xi(3t+1)\cdot e^{-3tx}+0\\
&-\frac{1}{3t+3}\xi(2)\cdot \xi(3t+1)\cdot e^{-3ty+x-y}\\
&-\frac{1}{3t}\frac{1}{3t+3}\cdot \xi(3t+2)\cdot e^{-3tx+x+2y}.\end{aligned}$$
Similarly, 
by taking residue along $z_2-z_3=1$ and assuming $z_3=s, z_2=s+1, z_1=-2s-1$, we get, using the tables in subsection B.1.2,
$$\begin{aligned}\omega_{\mathbb Q}^{SL(3)/P_{2,1};T}(t)=&-\frac{1}{3s+3}\xi(2)\cdot \xi(3s+1)\cdot e^{-3sx+x+2y}+0\\
&-\frac{1}{2}\frac{1}{3s+1}\cdot\xi(3s+3)\cdot e^{(3s+3)(x+y)}\\
&
+\frac{1}{2}\frac{1}{3s+2}\cdot\xi(3s+1)\cdot e^{-3sx}\\
&
-\frac{1}{3s}\frac{1}{3s+3}\cdot \xi(3s+2)\cdot e^{3sx+3sy+4x+2y}\\
&
+\frac{1}{3s}\xi(2)\cdot \xi(3s+3)\cdot e^{-3sy+x-y}.\end{aligned}$$
Clearly, there is no functional equation at this stage.
However, if we set $y=0$ in $T=(x,y,-x-y)$ so that $T=(x,0,-x)$, that is to say, $T=x\rho\in\mathbb C\cdot\rho$ sitting on the line spanned by $\rho$,
then we have
$$\boxed{\begin{aligned}\omega_{\mathbb Q}^{SL(3)/P_{1,2};x\rho}(t)=&\frac{1}{3t}\xi(2)\cdot \xi(3t+3)\cdot e^{3tx+4x}\\
&-\frac{1}{2}\frac{1}{3t+1}\cdot\xi(3t+3)\cdot e^{(3t+3)x}\\
&
+\frac{1}{2}\frac{1}{3t+2}\cdot\xi(3t+1)\cdot e^{-3tx}\\
&
-\frac{1}{3t+3}\xi(2)\cdot \xi(3t+1)\cdot e^{x}\\
&
-\frac{1}{3t}\frac{1}{3t+3}\cdot \xi(3t+2)\cdot e^{-3tx+x}\end{aligned}}$$
and
$$\boxed{\begin{aligned}\omega_{\mathbb Q}^{SL(3)/P_{2,1};x\rho}(t)=&-\frac{1}{3s+3}\xi(2)\cdot \xi(3s+1)\cdot e^{-3sx+x}\\
&
-\frac{1}{2}\frac{1}{3s+1}\cdot\xi(3s+3)\cdot e^{(3s+3)x}\\
&
+\frac{1}{2}\frac{1}{3s+2}\cdot\xi(3s+1)\cdot e^{-3sx}\\
&
-\frac{1}{3s}\frac{1}{3s+3}\cdot \xi(3s+2)\cdot e^{3sx+4x}\\
&
+\frac{1}{3s}\xi(2)\cdot \xi(3s+3)\cdot e^{x}\end{aligned}}$$
In particular, we have the functional equation
$$\boxed{\omega_{\mathbb Q}^{SL(3)/P_{1,2};x\rho}(-1-s)=\omega_{\mathbb Q}^{SL(3)/P_{2,1};x\rho}(x)}$$
Or put it in a better form, we set
\begin{equation}
\boxed{\begin{aligned}\xi_{\mathbb Q;{\bf T}}^{SL(3)/P_{1,2}}(s):=&\frac{1}{3s-3}\xi(2)\cdot \xi(3s)\cdot {\bf T}^{3s+1}-\frac{1}{3s}\xi(2)\cdot \xi(3s-2)\cdot {\bf T}\\
&
-\frac{1}{2}\frac{1}{3s-2}\cdot\xi(3s)\cdot {\bf T}^{3s}
+\frac{1}{2}\frac{1}{3s-1}\cdot\xi(3s-2)\cdot {\bf T}^{-3s+3}\\
&
-\frac{1}{3s-3}\frac{1}{3s}\cdot \xi(3t-1)\cdot {\bf T}^{-3s+4}\end{aligned}}
\end{equation}
and
\begin{equation}\boxed{\begin{aligned}\omega_{\mathbb Q;{\bf T}}^{SL(3)/P_{2,1}}(t)=&-\frac{1}{3s}\xi(2)\cdot \xi(3s-2)\cdot {\bf T}^{-3s+4}
+\frac{1}{3s-3}\xi(2)\cdot \xi(3s)\cdot {\bf T}\\
&
-\frac{1}{2}\frac{1}{3s-2}\cdot\xi(3s)\cdot {\bf T}^{3s}
+\frac{1}{2}\frac{1}{3s-1}\cdot\xi(3s-2)\cdot {\bf T}^{-3s+3}\\
&
-\frac{1}{3s-3}\frac{1}{3s}\cdot \xi(3s-1)\cdot {\bf T}^{3s+1}
\end{aligned}}
\end{equation}
Then we get
\begin{equation}
\boxed{\xi_{\mathbb Q;{\bf T}}^{SL(3)/P_{1,2}}(1-s)=\xi_{\mathbb Q;{\bf T}}^{SL(3)/P_{2,1}}(s)}
\end{equation}
This exposes a new symmetry for our zetas.
\vskip 0.5cm
\centerline{\bf REFERENCES}
\vskip 0.20cm
\noindent
[LS] J. Lagarias \& M. Suzuki, The Riemann Hypothesis for certain 
integrals of Eisenstein series, {\it J Number Theory}, 118(2006) 98-122
\vskip 0.15cm
\noindent
[S] M. Suzuki, A proof of the Riemann Hypothesis for the Weng zeta function of rank 3 for the rationals,  in {\it Conference on L-Functions}, 175-200, World Sci. (2007)
\vskip 0.15cm
\noindent
[S2]  M. Suzuki, The Riemann hypothesis for Weng's zeta function of ${\rm Sp}(4)$ over $\mathbb Q$, preprint, 2008,
available at http://xxx.lanl.gov/abs/0802.0102
\vskip 0.15cm
\noindent
[SW] M. Suzuki \& L. Weng, Zeta functions for $G_2$ and their zeros,  preprint, 2008, available at http://xxx.lanl.gov/abs/0802.0104
\vskip 0.15cm
\noindent
[W0] L. Weng, Non-abelian zeta function for function fields, {\it Amer. J. Math} 127 (2005), 973-1017
\vskip 0.15cm
\noindent
[W1] L. Weng, Geometric Arithmetic: A Program, in {\it Arithmetic Geometry and Number Theory},  211-390, World Sci. (2006)
\vskip 0.15cm
\noindent
[W2] L. Weng,  A Rank two zeta and its zeros, {\it J of Ramanujan Math. Soc.}, 21 (2006), 205-266
\vskip 0.15cm
\noindent
[W3]  L. Weng, A geometric approach to $L$-functions, in {\it Conference on L-Functions}, 219-370, World Sci (2007)
\vskip 0.15cm
\noindent  
{[W4]} L. Weng, Zeta function for $Sp(2n)$, 
Appendix to [S2],  preprint, 2008, available at http://xxx.lanl.gov/abs/0802.0102
\eject
\vskip 0.5cm
\noindent
{\bf Lin WENG}\footnote{
{\it Acknowledgement.}  Special thanks due to Deninger and Ueno for their constant encouragements, due to Henry Kim for bringing to our attention the paper of Diehl, and due to Suzuki for testing numerically the RH associated to $SL(4,5),$  $Sp(4)$ and $G_2$.

This work is partially supported by JSPS.} 

 
\smallskip
\noindent
Graduate School of Mathematics

\noindent
Kyushu University

\noindent 
Fukuoka 812-8581

\noindent
Japan

\noindent
Email: weng@math.kyushu-u.ac.jp

\smallskip
\noindent
{\it and}

\smallskip
\noindent
Chennai Mathematical Institute

\noindent
Plot H1, SIPCOT IT Park

\noindent
Padur PO, Siruseri 603103

\noindent
 India 
\end{document}